\documentclass[oupthm,crop,info]{CUP-JNL-CJM}

\usepackage{amsmath,amssymb,amsfonts}
\usepackage{bbm,bm}
\usepackage{tikz-cd}
\usepackage{mathrsfs}

\usepackage{float}
\usepackage{graphicx}
\usepackage{accents}
\usepackage{rotating}
\usepackage{appendix}
\usepackage[numbers]{natbib}
\usepackage{graphicx}

\usepackage{booktabs,extarrows,enumerate,mathtools}

\theoremstyle{oupplain}
\newtheorem{theorem}{Theorem}[section]
\newtheorem{lemma}[theorem]{Lemma}
\newtheorem{corollary}[theorem]{Corollary}
\newtheorem{proposition}[theorem]{Proposition}
\theoremstyle{oupdefinition}

\theoremstyle{oupremark}
\newtheorem{remark}[theorem]{Remark}
\newtheorem{example}[theorem]{Example}
\theoremstyle{oupproof}
\newtheorem{proof}{Proof}

\numberwithin{equation}{section}

\newcommand{\loca}{\ensuremath{(\frac{1}{2})}}
\newcommand{\Z}{\ensuremath{\mathbb{Z}}}
\newcommand{\ZZ}[1]{\ensuremath{\Z_{(#1)}}}
\newcommand{\z}[1]{\ensuremath{\Z/2^{#1}}}
\newcommand{\zp}[1]{\ensuremath{\Z/p^{#1}}}
\newcommand{\xra}{\xlongrightarrow}
\newcommand{\lra}[1]{\ensuremath{\langle #1 \rangle}}
\newcommand{\PP}{\ensuremath{\mathcal{P}^1}}

\DeclareMathOperator{\Sq}{Sq}
\DeclareMathOperator{\im}{im}
\DeclareMathOperator{\coker}{coker}
\DeclareMathOperator{\Hom}{Hom}
\DeclareMathOperator{\sk}{sk}

\newcommand{\vees}[1]{\ensuremath{\mathop{\bigvee}\limits_}}

\newcommand{\mat}[4]{\ensuremath{\begin{bmatrix}
  #1&#2\\
  #3&#4
\end{bmatrix}}}

\newcommand{\matwo}[2]{\ensuremath{{\begin{bmatrix}
  #1\\
  #2
\end{bmatrix}}}}

\articletype{RESEARCH ARTICLE}
\jname{Canadian Mathematical Society}
\jyear{2020}

\begin{document}

\begin{Frontmatter}

\title{Suspension Splitting and Cohomotopy Sets of Simply Connected $7$-manifolds} 

\author[1]{Ruizhi Huang} 
\author[2]{Pengcheng Li}
\authormark{R. Huang and P. Li}

\address[1]{\orgname{Institute of Mathematics, Academy of Mathematics and Systems Science, 
 Chinese Academy of Sciences}, \orgaddress{\city{Beijing}, \country{China}} 
\email{huangrz@amss.ac.cn}} 

\address[2]{\orgname{Department of Mathematics, School of Sciences, Great Bay University}, \orgaddress{\city{Dongguan}, \country{China}} 
\email{lipcaty@outlook.com}}



\abstract{
Let $M$ be a closed simply connected $7$-manifold. In this paper we establish homotopy decompositions of the reduced suspension space $\Sigma M$ into a wedge sum of simpler spaces when localized at a set of primes. These decompositions are applied to study the cohomotopy sets $\pi^k(M)$ and the $p$-local cohomotopy sets $\pi^4(M;\mathbb{Z}_{(p)})$. }

\end{Frontmatter}

\section{Introduction}\label{sec:intro}

Seven-dimensional manifolds play a central role in geometry and topology. In differential geometry, for example, there has been growing interest in $G_2$-manifolds, particularly following the seminal work of Donaldson and Thomas \cite{DT}.
Changing the focus to the algebraic topology of gauge theory, a general principle \cite{The, So, Hua2} posits that the homotopy theory of gauge groups is closely related to the homotopy type of the underlying manifolds.  In geometric topology, the foundational contributions by Kreck \cite{Kre2} and Crowley-Nordstr\"{o}m \cite{CN} independently established a complete classification of closed smooth $2$-connected $7$-manifolds using systems of algebraic invariants. Building on his earlier breakthrough work \cite{Kre1}, Kreck \cite{Kre2} further extended this classification to simply connected 7-manifolds with torsion-free second homology.

In recent years, the homotopy theory of suspended manifolds has attracted significant interest among topologists. For example, So and Theriault \cite{ST} and Li \cite{lpc4} determined the homotopy type of suspended $4$-manifolds. Subsequent work by Huang \cite{Hua1}, Li and Zhu \cite{LZ}, and Amelotte, Cutler, and So \cite{ACS} extended this analysis to certain $5$-manifolds. Huang \cite{Hua3} and Cutler and So \cite{CS} further advanced these results for $6$-manifolds. Additionally, Membrillo-Solis \cite{MS} investigated the problem in the context of $S^3$-bundles over $S^4$ at large primes. Collectively, these studies have applications in characterizing key invariants in geometry and mathematical physics, such as reduced $K$-groups and gauge groups.

This paper focuses on the homotopy type of the (reduced) suspension space of $7$-manifolds.
Let $M$ be a closed simply connected $7$-manifold whose reduced integral homology groups $H_\ast(M;\Z)$ are given by the following table:  
\begin{equation}\label{HMeq}
  \begin{tabular}{ccccccccc}
    \toprule
$i$& $2$&$3$&$4$&$5$&$0,7$&$\text{otherwise}$
  \\\midrule
  $H_i(M;\Z)$& $\Z^l \oplus H$&  $\Z^k\oplus T$ & $\Z^k \oplus H $& $\Z^l$&$\Z$&0
  \\ \bottomrule
  \end{tabular}~~,
\end{equation}
where $k,l\geq 0$, and $H$ and $T$ are two finitely generated abelian torsion groups.
For integer $n\geq 2$ and finitely generated abelian group $A$, let $P^n(A)$ be the Moore space such that the reduced cohomology $\widetilde{H}^i(P^n(A);\Z)$ is isomorphic to the abelian group $A$ if $i=n$ and is trivial otherwise \cite{N2}. For each $k>1$, denote $P^n(k)=P^n(\Z/k)$. There is a homotopy cofibration sequence:
\[S^{n-1}\xra{k}S^{n-1}\xra{i_{n-1}}P^n(k)\xra{q_n}S^n,\]
where $k$ is the degree $k$ map, $i_{n-1}$ and $q_n$ are the canonical inclusion and pinch maps, respectively. 
For based CW-complexes $X,Y$ and each prime $p$, we denote by $X_{(p)}$ the $p$-localized space of $X$, denote by $X\simeq_{(p)} Y$ if $X$ is homotopy equivalent to $Y$ when localized at $p$, denote by $X\simeq_{\loca}Y$ if $X$ is homotopy equivalent to $Y$ when localized away from $2$. By Lemma \ref{lem:pi6p3}, for each $r$ there is a lift $\widetilde{\alpha}_r\colon S^6\to P^3(3^r)$ such that $q_3\circ \widetilde{\alpha}_r=\alpha$, where $\alpha\in \pi_6(S^3)$ is a generator of order $3$. Recall the Steenrod reduced power operation $\PP=\PP_3\colon H^\ast(-;\Z/3)\to H^{\ast+4}(-;\Z/3)$.
By naturality, the fact that $\alpha$ is detected by $\PP$ \cite{Har} implies that the maps $\widetilde{\alpha}_{r}$, $i_3\alpha$ and their suspensions are also detected by $\PP$; see Lemma \ref{lem:PP} for more details. For each prime $p$ and each integer $r\geq 1$, the $r$-th higher Bockstein operation \[\beta_r\colon H^\ast(-;\Z/p)\dashrightarrow H^{\ast+1}(-;\Z/p)\] is inductively defined by setting $\beta_1$ as the usual Bockstein homomorphism associated to the short exact
sequence 
$0\to \Z/p\to \Z/p^2\to \Z/p\to 0;$
for $r\geq 2$, $\beta_r$ is defined on the intersection of $\ker(\beta_i)$,  $i<r$, and takes values in the quotient by the $\im (\beta_i)$, $i<r$; this is also indicated by the dashed arrow above. See \cite[Section 5.2]{Har} for more details. 

Our first main result classifies the homotopy type of the suspension space $\Sigma M$ after localized away from $2$. Note that we allow that $k=l=0$ in (\ref{HMeq}), and we adopt the convention that $\bigvee_{i=1}^{t}X_i=\ast$ if $t<1$ in Theorem \ref{thm:main}.

\begin{theorem}\label{thm:main}
  Let $M$ be a closed simply connected $7$-manifold with $H_\ast(M;\Z)$ given by (\ref{HMeq}). 
  \begin{enumerate}[1.]
    \item If $\PP\colon H^3(M;\Z/3)\to H^7(M;\Z/3)$ is trivial, then there is a homotopy equivalence 
     \[ \Sigma M\simeq_{\loca}  \bigvee_{i=1}^l(S^3\vee S^6)\vee \bigvee_{j=1}^k (S^4\vee S^5)  \vee  P^4(H) \vee P^5(T)\vee P^6(H)\vee S^8.\]
    \item If $\PP\colon H^3(M;\Z/3)\to H^7(M;\Z/3)$ is non-trivial, then the homotopy type of $\Sigma M$ can be described as follows.
    \begin{enumerate} 
    \item\label{SM:PP=1:P4} If there exists $x\in H^3(M;\Z/3)$ such that $\PP(x)\neq 0$ and $x\in \im (\beta_{r})$ for some $r\geq 1$, then there is a homotopy equivalence 
    \[\Sigma M\simeq_{\loca}\bigvee_{i=1}^l(S^3\vee S^6)\vee \bigvee_{j=1}^k (S^4\vee S^5)  \vee  P^4\big(\frac{H}{\Z/3^{\bar{r}}}\big) \vee P^5(T) \vee P^6(H)\vee  X^8(\widetilde{\alpha}_{\bar{r}}),\]
    where $ X^8(\widetilde{\alpha}_{\bar{r}})=P^4(3^{\bar{r}})\cup_{\Sigma\widetilde{\alpha}_{\bar{r}}}e^8$ and $\bar{r}$ is the minimum of $r$ such that $x\in \im(\beta_r)$ and $\PP(x)\neq 0$.
    \item\label{SM:PP=1:S4} If for any $x\in H^3(M;\Z/3)$ with $\PP(x)\neq 0$, there hold $x\notin \im(\beta_r)$ and $\beta_s(x)=0$ for any $r,s\geq 1$, then 
    there is a homotopy equivalence 
    \[\Sigma M\!\simeq_{\loca}\!\bigvee_{i=1}^l\!(S^3\vee S^6)\!\vee \!\bigvee_{j=1}^{k-1}\! S^4\!\vee \!\bigvee_{j'=1}^{k}\!S^5 \! \vee\!  P^4(H) \!\vee\! P^5(T) \!\vee\! P^6(H)\!\vee\! (S^4\cup_{\Sigma\alpha}e^8).\]
    \item\label{SM:PP=1:P5} If for any $x\in H^3(M;\Z/3)$ with $\PP(x)\neq 0$, there hold  $x\notin \im(\beta_r)$ for any $r\geq 1$, while $\beta_s(x)\neq 0$ for some $s\geq 1$, then there is a homotopy equivalence 
    \[\Sigma M\simeq_{\loca}\bigvee_{i=1}^l(S^3\vee S^6)\vee \bigvee_{j=1}^k (S^4\vee S^5)  \vee  P^4(H)  \vee P^5(\frac{T}{\Z/3^{\bar{s}}})\vee P^6(H) \vee X^8_{\bar{s}},\]
    where $X^8_{\bar{s}}=P^5(3^{\bar{s}})\cup_{i_4(\Sigma\alpha)}e^8$ and $\bar{s}$ is the maximum of $s$ such that $\beta_s(x)\neq 0$.
  \end{enumerate}
  \end{enumerate}
\end{theorem}

Recall that there is a secondary operation $\Theta$ based on the relation 
\[\Sq^1(\Sq^2\Sq^1)+\Sq^2\Sq^2=0\]
that detects $\eta^2\in \pi_{n+2}(S^n)$; see \cite[Page 96]{Har}.
The $2$-local homotopy type of the suspension $\Sigma M$ is very difficult to handle, we have the following description of the double suspension $\Sigma^2 M_{(2)}$ in a special case. 

\begin{theorem}[see Theorem \ref{thm:SM-2-local}]\label{thm:2-local}
	Let $M$ be a closed simply connected $7$-manifold with $H_\ast(M;\Z)$ given by (\ref{HMeq}), where $H$ and $T$ are both $2$-torsion-free. Suppose that $\Sq^2$ acts trivially on $H^2(M;\z{})$ and $H^2(M;\z{})$ and that $\Theta$ acts trivially on $H^\ast(M;\z{})$, then there is a homotopy equivalence 
\[\Sigma^2 M\simeq_{(2)}\bigvee_{r=1}^k S^6\vee \bigvee_{s=1}^{l-b}S^7\vee \bigvee_{i=1}^{l-1}S^4\vee \bigvee_{j=1}^{k-b-1}S^5\vee \bigvee_{v=1}^{b-1}C^7_\eta\vee C_{\hbar},\]
where $C^7_\eta=\Sigma^3\mathbb{C}P^2$, $C_{\hbar}$ is the homotopy cofibre of some map $\hbar\colon S^8\to S^4\vee S^5\vee C^7_\eta$ given by (\ref{eq:hbar}) and $0\leq b\leq \min\{k,l\}$ is an integer depending on $M$; moreover, $b=0$ if and only if $\Sq^2$ acts trivially on $H^3(M;\z{})$.
\end{theorem}

As demonstrated in \cite{ST,Hua1,CS}, the homotopy type of the suspension space determines the homotopy type of gauge groups of certain principal bundles over $M$. Theorems \ref{thm:main} and \ref{thm:2-local} have potential applications in characterizing the ($p$-local) homotopy types of the gauge groups of simply connected $7$-manifolds. 

By the Pontryagin-Thom construction (compare \cite[Chapter IX, Theorem 5.5]{Kosinski93}) , the characterization of the cohomotopy set $\pi^k(M)=[M,S^k]$,  which consists of homotopy classes of based maps from a smooth closed manifold $M$ to the $k$-sphere, has become an actively studied topic in geometry and topology,  see \cite{Taylor12,KMT12} for a complete characterization of $\pi^k(M)$ of $4$-manifolds $M$. 

The homotopy type of $\Sigma M$ is closely related to the cohomotopy sets $\pi^k(M)$ and the $p$-local cohomotopy sets $\pi^k(M;\ZZ{p})=[M,S^{k}_{(p)}]$, where $p$ is a prime. On the one hand, the suspension isomorphism $\pi^k(M)\cong \pi^{k+1}(\Sigma M)$ holds in the stable range $\dim(M)\leq 2k-2$. On the other hand, unstable homotopy theory  emerges as a powerful tool for characterizing the cohomotopy sets $\pi^k(M)$ outside the stable range (see, for example, \cite{LPW,LZ,ACS}). 

For a closed simply connected $7$-manifold $M$, we say that $M$ is \emph{spin} if the reduced Steenrod square $\Sq^2$ acts trivially on $H^5(M;\z{})$; otherwise we say that $M$ is \emph{nonspin}. When the manifold $M$ is smooth, being spinnable is equivalent to the second Stiefel-Whitney class $w_2(M)$ vanishes, by the Wu's formula $\Sq^2(x)=x\smallsmile w_2(M)$ \cite{Wu} for any $x\in H^5(M;\z{})$. In the second part of this paper, we apply Theorem \ref{thm:main} to study the cohomotopy sets $\pi^\ast(M)$ of a closed simply connected $7$-manifold $M$. Note that $\pi^k(M)$ is an abelian group in the stable range $k\geq 5$, by the generalized Freudenthal's suspension theorem (compare \cite[Theorem 1.21]{Cohen70}).
For each $k\geq 1$, there are the  ($p$-local) cohomotopy Hurewicz maps 
\[h^k\colon \pi^k(M)\to H^k(M;\Z),\quad  h^k_{(p)}\colon \pi^k(M;\ZZ{p})\to H^k(M;\ZZ{p})\]
induced by the canonical inclusions $S^k\hookrightarrow K(\Z,k)$ and $S^k_{(p)}\hookrightarrow K(\Z_{(p)},k)$, respectively. Moreover, the Hurewicz maps $h^k$ and $h^k_{(p)}$ are homomorphisms for $k\geq 5$. The following theorem summarizes the main results of Section \ref{sec:chtp}.

\begin{theorem}\label{thm:chtp}
	Let $M$ be a closed simply connected $7$-manifold with $H_\ast(M;\Z)$ given by (\ref{HMeq}). Set $\varepsilon=0$ if $M$ is spin, otherwise $\varepsilon=1$.
	\begin{enumerate}[1.]
    \item\label{chtp:dim=1} $\pi^k(M)=0$  for $k=1$ or $k\geq 8$.
		\item\label{chtp:dim=6}  There are group isomorphisms $\pi^7(M)\cong\Z$ and $\pi^6(M)\cong\z{1-\varepsilon}$.
		\item\label{chtp:dim=5}  If the manifold $M$ is smooth and the torsion group $H$ is $2$-torsion-free, then there is a split exact sequence of groups
    \[0\to \z{1-\varepsilon}\to \pi^5(M)\xra{\jmath}\ker(\Sq^2_\Z)\to 0,\]
    where $\ker(\Sq^2_\Z)\cong H^5(M;\Z)$ and $\jmath$ is the homomorphism such that the composition 
    \[\pi^5(M)\xra{\jmath}\ker(\Sq^2_\Z)\cong H^5(M;\Z)\]
    is the fifth cohomotopy Hurewicz homomorphism. Here $\Sq^2_\Z=\Sq^2\circ \rho_2$ is composition of the reduced Steenrod square $\Sq^2$ with the mod $2$ reduction $\rho_2\colon H^5(M;\Z)\to H^5(M;\z{})$.
    
    \item\label{chtp:dim=3}  The suspension $E_\ast\colon \pi^3(M)\to \pi^4(\Sigma M)$ is a group isomorphism in the following two cases:
    \begin{enumerate}
      \item the manifold $M$ is nonspin;
      \item all groups are localized away from $2$, that is, $E_\ast\colon \pi^3(M)_{\loca}\xra{\cong} \pi^4(\Sigma M)_{\loca}$.
    \end{enumerate}
   Moreover, if the conditions in Theorem \ref{thm:main} (\ref{SM:PP=1:P4}) holds, then there holds a group isomorphism
   \[\pi^3(M)\cong_{\loca}
    \Z^k\oplus (H/(\Z/3^{\bar{r}}))\oplus \Z/3^{\bar{r}-1},\]
    where $\bar{r}$ is given by Theorem  \ref{thm:main} (\ref{SM:PP=1:P4}); otherwise there is a group isomorphism 
    \[\pi^3(M)\cong_{\loca}H^3(M;\Z)\cong_{\loca} \Z^k\oplus H.\]
    \item\label{chtp:dim=4} For each prime $p\geq 5$, the $p$-local cohomotopy Hurewicz map 
    \[h_{(p)}^4\colon \pi^4(M;\ZZ{p})\to H^4(M;\ZZ{p})\] 
    is surjective.  For $u\in H^4(M;\ZZ{p})$, let $e\in \pi^4(M;\ZZ{p})$ be such that $h_{(p)}^4(e)=u$.  There is a bijection between $(h_{(p)}^4)^{-1}(u)$ and the cokernel of the homomorphism 
    \[\psi_e\colon H^3(M;\ZZ{p})\to H^7(M;\ZZ{p}),\quad \psi_e(v)=v\smallsmile e^\ast(\iota),\] 
    where $\iota\in H^4(S^4_{(p)};\ZZ{p})$ is the fundamental class.  
		\end{enumerate}
\end{theorem}

When $M$ is smooth and spin, that is, $w_2(M)=0$, the suspension $E_\ast$ in Theorem \ref{thm:chtp} (\ref{chtp:dim=3}) is usually not surjective, see Example  \ref{ex:notsurj}. Theorem \ref{thm:chtp} (\ref{chtp:dim=4}) is also generalized to characterize the $p$-local cohomotopy set $\pi^4(X;\ZZ{p})$ of a CW-complex $X$ of dimension at most $10$, see Theorem \ref{thm:chtp=4}.

The paper is organized as follows.  Section \ref{sec:prelim}  is divided into two subsections to give technical computations of homotopy groups of mod $p^r$ Moore spaces and suspensions of $\mathbb{C}P^2$, review useful criteria for null homotopy, respectively. Section \ref{sec: hdec} focuses on the application of homology decomposition techniques to analyze the reduced suspension space of a simply connected $7$-manifold $M$, assuming that $M$ contains no $2$-torsion in its homology. Moving forward, Section \ref{sec:SM:p-local} delves into the homotopy type of the suspension $\Sigma M$ when localized away from $2$, culminating in the proof of Theorem \ref{thm:main}. At the end of this section, we study the homotopy type of the double suspension space $\Sigma^2M$ localized at $2$ and prove Theorem \ref{thm:SM-2-local}. Section \ref{sec:chtp} is divided into three subsections to investigate ($p$-local) cohomotopy sets for a simply connected $7$-manifold, leading to the proof of Theorem \ref{thm:chtp}.

\medskip

\noindent{\bf Acknowledgments.} 
The authors would like to thank the anonymous referee for his careful reading and valuable comments that improved the quality of this paper. 
Ruizhi Huang was supported in part by the National Natural Science Foundation of China (Grant nos. 12331003 and 12288201), the National Key R\&D Program of China (No. 2021YFA1002300), the Youth Innovation Promotion Association of Chinese Academy Sciences, and the ``Chen Jingrun'' Future Star Program of AMSS.
Pencheng Li was supported by the National Natural Science Foundation of China (Grant No. 12101290).
 

\section{Preliminaries}
\label{sec:prelim}

Throughout the paper, all spaces are based CW-complexes, all maps are base-point-preserving and are identified with their homotopy classes in notation. 
For an abelian group $G$, denote by  
\[G\cong C_1\langle x_1\rangle\oplus\cdots\oplus C_n\langle x_n\rangle\]  
if $G$ has generators $x_1,\cdots,x_n$, where $C_i$ are non-trivial cyclic groups. For integers $a$ and $b$, we denote by $(a,b)$ the greatest common divisor of $a$ and $b$.
We shall frequently use the following homotopy groups and notation (compare \cite{Tod}):
\begin{itemize}
  \item $\pi_3(S^2)\cong\Z\lra{\eta_2}$, $\pi_{n+1}(S^n)\cong \z{}\lra{\eta_n}$ for $n\geq 3$, where $\eta_n=\Sigma^{n-2}\eta_2$ for $n\geq 2$; 
  \item  $\pi_{n+2}(S^n)\cong\z{}\lra{\eta^2}$, where $\eta^2=\eta_{n}\eta_{n+1}$; 
  \item $\pi_6(S^3)\cong\Z/12\lra{\nu'}$, $\pi_7(S^4)\cong \Z\lra{\nu_4}\oplus\Z/12\lra{\Sigma \nu'}$, and $\pi_{n+3}(S^n)\cong\Z/24\lra{\nu_n}$ for $n\geq 5$, where $\nu_{n}=\Sigma^{n-4}\nu_4$ for $n\geq 4$. We have $\eta^3=6\Sigma\nu'$, $\Sigma^2\nu'=2\nu_5$.
\end{itemize}
Note that if there is no confusion, we use the same notation $\eta$ to denote the iterated suspension $\eta_n$ for different $n$.

\subsection{Homotopy groups of mod $p^r$ Moore spaces and $C^{n+2}_\eta$}
For integers $n\geq 3$ and $r,s\geq 1$, we denote by 
\[i_n^r\colon S^n\to P^{n+1}(p^r) \text{ and }q_{n+1}^s\colon P^{n+1}(p^s)\to S^{n+1}\]  the canonical inclusion map of the bottom cell into $P^{n+1}(p^r)$ and the canonical pinch map from $P^{n+1}(p^s)$ onto $S^{n+1}$ by collapsing the bottom cell, respectively. If there is no confusion, we will omit the superscripts $r$ and $s$ in the notation. 
We have the following classification of homotopy classes of based maps between Moore spaces and spheres.
\begin{lemma}[compare \cite{BH91}]\label{lem:Moore1}
  Let $p$ be an odd prime and let $m,n\geq 3$ and $r,s\geq 1$ be integers. 
\begin{enumerate}
  \itemsep=2pt
  \item\label{SkPr} $[S^n,P^{n+1}(p^r)]\cong\zp{r}\langle i_n^r\rangle$ and $[P^{n+1}(p^s),S^{n+1}]\cong\zp{s}\langle q_{n+1}^s\rangle$.
  \item $[S^{n+1},P^{n+i}(p^r)]=[P^{n+1+i}(p^r),S^{n}]=0$ for $i=0,1$. 
 \item $[P^{n+1}(p^r),P^n(p^s)]=0$ and $[P^n(p^r),P^m(q^s)]=0$ for coprime $p,q$.
  \item\label{PrPs} $[P^{n+1}(p^r),P^{n+1}(p^s)]\cong \zp{\min\{r,s\}}\lra{B_n(\chi^r_s)}$, where $B_n(\chi^r_s)$ satisfies the formulae
  \begin{equation}\label{eq:chi}
    \begin{aligned}
      B_n(\chi^r_s)i_{n}^r\simeq\left\{\begin{array}{ll}
    i_{n}^s,&r\geq s;\\
    p^{s-r}i_{n}^s,&r\leq s.
  \end{array}\right.& ~q_{n+1}^sB_n(\chi^r_s)\simeq\left\{\begin{array}{ll}
    p^{r-s}q_{n+1}^r,&r\geq s;\\
    q_{n+1}^r,&r\leq s.
  \end{array}\right.
    \end{aligned}
  \end{equation}

\end{enumerate}
\end{lemma}
Note that $B_n(\chi^r_s)$ satisfies the formula $\Sigma B_n(\chi^r_s)=B_{n+1}(\chi^r_s)$ and induces the homomorphism of the $n$-th integral homology groups 
\begin{equation}\label{eq:chi^r_s}
  \chi^r_s\colon \Z/p^r\to\Z/p^s
\end{equation}
which is the reduction map if $r\geq s$ and the multiplication by $p^{s-r}$ if $r<s$.


Let $S^n\{p^r\}$ be the homotopy fibre of the degree map $p^r$ on $S^n$ and let $F^{2n+1}\{p^r\}$ be the homotopy fibre of the pinch map $q_{2n+1}\colon  P^{2n+1}(p^r)\to S^{2n+1}$.

\begin{lemma}\label{Lem:CMN}
Let $p$ be an odd prime and let $n\geq 1$.
\begin{enumerate}
  \item\label{nei87} There is a $p$-local homotopy equivalence 
  \begin{equation}\label{cmn-odd}
    \Phi\colon T^{2n+1}\{p^r\}\times \Omega\Sigma \big(\bigvee_\alpha P^{n_\alpha}(p^r)\big)\xra{\simeq_{(p)}}\Omega P^{2n+1}(p^r),
  \end{equation}
  where $n_\alpha\geq 4n-1$ and for each $n$, there exists exactly one $\alpha$ such that $n_\alpha=4n-1$; the space $T^{2n+1}\{p^r\}$ sits in the $p$-local homotopy fibration sequences
  \begin{equation}\label{fib:T}
    S^{2n-1}\times \prod_{k=1}^{\infty}S^{2p^kn-1}\{p^{r+1}\} \xra{\partial_r} T^{2n+1}\{p^r\}\to \Omega S^{2n+1},
  \end{equation}
\begin{equation}\label{fib:TT}
  W_n\times \prod_{k=1}^{\infty}S^{2p^kn-1}\{p^{r+1}\} \to T^{2n+1}\{p^r\}\to \Omega S^{2n+1}\{p^r\},
\end{equation}
where $W_n$ is the homotopy fibre of the double suspension $E^2\colon S^{2n-1}\to \Omega^2S^{2n+1}$, and $W_n$ is $(2pn-4)$-connected.
  \item\label{nei81} There is a $p$-local homotopy equivalence 
  \begin{equation}\label{CMN:OmegaF}
    \Psi\colon S^{2n-1}\times \prod_{k=1}^{\infty} S^{2p^kn-1}\{p^{r+1}\}\times \Omega\Sigma\big(\bigvee_{\alpha} P^{n_\alpha}(p^r)\big)\xra{\simeq_{(p)}}\Omega F^{2n+1}\{p^r\},                
  \end{equation}
where $\bigvee_{\alpha} P^{n_\alpha}(p^r)$ is an infinite bouquet of mod $p^r$ Moore spaces, with only finitely many Moore spaces in each dimension and the least value of $n_\alpha$ being $4n-1$. 
\end{enumerate} 
\begin{proof}
 The homotopy equivalence (\ref{cmn-odd}) refers to \cite[Corollary 1.9]{CMN87} or \cite[Proposition 0.3 (b)]{Neisen87},  the homotopy fibre sequences (\ref{fib:T}) and (\ref{fib:TT}) refer to \cite[the main diagram on Page 38]{Neisen87}.
The homotopy equivalence (\ref{CMN:OmegaF}) refers to \cite[Theorem 12.1]{CMN79} or \cite[Theorem 3.1]{Neisen81}.
\end{proof}
\end{lemma}

\begin{lemma}\label{lem:pi7P5}
Let $p$ be an odd prime and let $r$ be a positive integer. There is a group isomorphism 
 \[\pi_7(P^5(p^r))\cong \zp{r}\lra{[1_P,1_P]\circ \widetilde{i_7}}\oplus \Z/(p,3)\lra{i_4(\Sigma\alpha)},\] 
 where $[1_P,1_P]$ denotes the Whitehead product of the identity $1_P$ on $P^4(p^r)$ and $\widetilde{i_7}\in \pi_7(\Sigma P^4(p^r)\wedge P^4(p^r))$ is the generator corresponding to the generator $i_7\in \pi_7(P^8(p^r))$.

  \begin{proof} 
   There are group isomorphisms by taking $n=2$ in (\ref{cmn-odd}): 
    \begin{align*}
      \pi_7(P^5(p^r))\cong &\pi_6(\Omega P^5(p^r))\cong \pi_6(T^5\{p^r\})\oplus \zp{r}.
    \end{align*}
    Since $S^{2n+1}\{p^r\}$ is $(2n-1)$-connected and $W_n$ is $(2pn-4)$-connected, we have 
    \[\pi_i(W_n\times \prod_{k=1}^{\infty}S^{2p^kn-1}\{p^{r+1}\})=0 \text{ for }i=5,6 \text{ and }n=2.\]
    Then by the homotopy fibration sequence (\ref{fib:TT}), we get 
    \[\pi_6(T^5\{p^r\})\cong \pi_6(\Omega S^5\{p^r\})\cong \Z/(3,p).\]
    Thus $\pi_7(P^5(p^r))\cong\zp{r}\oplus\Z/(3,p)$.

Let $\Phi_1$ and $\Phi_2$ be the restrictions of the homotopy equivalence $\Phi$ in (\ref{cmn-odd}) to $T^5\{p^r\}$ and $\Omega P^8(p^r)$, respectively. By the construction of $\Phi$, there is a homotopy commutative diagram 
    \[\begin{tikzcd}
    \pi_6(P^4(p^r)\wedge P^4(p^r))\ar[rr,"{\langle E,E\rangle_\ast}"]\ar[d,"E_\ast","\cong"swap]&& \pi_6(\Omega\Sigma P^4(p^r))\\
    \pi_6(\Omega\Sigma P^4(p^r)\wedge P^4(p^r))\ar[urr,"(\Phi_2)_\ast"swap]&&
    \end{tikzcd},\]
where $\lra{E,E}$ denotes the Samelson product of the suspension $E\colon P^4(p^r)\to \Omega \Sigma P^4(p^r)$, and the vertical $E_\ast$ is an isomorphism by the Freudenthal suspension theorem. Since the adjoint of $\lra{E,E}$ is the Whitehead product $[1_P,1_P]$, we get $(\Phi_2)_\ast(\widetilde{i_7})= [1_P,1_P]\circ \widetilde{i_7}$, which generates the direct summand $\zp{r}$ of $\pi_7(P^5(p^r))$.

When $p=3$, the homotopy fibration sequence (\ref{fib:T}) with $n=2, p=3$ induces an exact sequence of groups localized at $3$:
    \[\pi_7(\Omega S^5)\xra{}\pi_6(S^3)\xra{(\partial_r')_\ast}\pi_6(T^5\{3^r\})\to 0,\]
where the first unlabelled map is the connecting homomorphism and  $\partial_r'$ is the composition \[\partial_r'\colon S^3\hookrightarrow S^3\times \prod_{i=1}^{\infty}S^{2\cdot 3^kn-1}\{3^{r+1}\}\xra{\partial_r}T^5\{3^r\}.\] 
Since $\pi_6(S^3)_{(3)}\cong\Z/3$ and $\pi_6(T^5\{3^r\})\cong \Z/3$, the induced homomorphism $(\partial_r')_\ast$ is an isomorphism. 
On the other hand, the composition 
\[\Phi_1'\colon S^3\xra{\partial_r'}T^5\{3^r\}\xra{\Phi_1}\Omega P^5(3^r)\]
is homotopic to the inclusion of the bottom cell; hence we have 
\[\Phi_1'\simeq (\Omega i_4)\circ E\colon S^3\xra{E}\Omega S^4\xra{\Omega i_4}\Omega P^5(3^r).\]
It follows that there is a split monomorphism 
\[\pi_6(S^3)_{(3)}\cong\Z/3\lra{\alpha}\xra{E_\ast}\pi_7(S^4)_{(3)}\xra{(i_4)_\ast}\pi_7(P^5(3^r)).\]
Thus we get the group $\pi_7(P^5(p^r))$ and its generators as stated in the Lemma. 
  \end{proof}
\end{lemma}

\begin{lemma}\label{pin+2pnlemma}
For any $n\geq 6$, we have 
\[
\pi_{n+2}(P^n(p^r))=0~{\rm for}~p\geq 5, \quad \pi_{n+2}(P^n(3^r))\cong\Z/3 \lra{ i_{n-1}(\Sigma^{n-4}\alpha)}.
\]
  
\end{lemma}
\begin{proof}
The statement for $p\geq 5$ is proved by \cite[Lemma 6.5]{Hua3}. Let $p=3$.
By the Freudenthal suspension theorem, it suffices to show $\pi_{9}(P^7(3^r))\cong\Z/3\lra{ i_6(\Sigma^{3}\alpha)}$. 

Consider the cofibration $S^6\stackrel{ 3^r}{\longrightarrow}S^6\xra{i_6}P^7(3^r)$. By the Blakers-Massey theorem (compare \cite{BM}), it is a homotopy fibration up to degree $10$. Hence, there is an exact sequence
\[
\Z/3\cong\pi_9(S^6)_{(3)}\stackrel{ 3^r}{\longrightarrow}\pi_9(S^6)_{(3)}\xra{i_{6\ast}}\pi_{9}(P^7(3^r))_{(3)}\to \pi_8(S^6)_{(3)}=0,
\]
where the multiplication map $3^r$ is trivial. Hence $\pi_{9}(P^7(3^r))_{(3)} \cong \pi_9(S^6)_{(3)}\cong\Z/3\lra{\Sigma^3 \alpha}$. Since $P^7(3^r)$ is contractible at any prime $p\neq 3$ and rationally, $\pi_{9}(P^7(3^r)) \cong\pi_{9}(P^7(3^r))_{(3)} \cong\Z/3\lra{ i_6(\Sigma^{3}\alpha)}$ as required.
This completes the proof of the lemma.
\end{proof}

\begin{lemma}
  \label{lem:pi6p3}
  There is a group isomorphism 
  \[\pi_6(P^3(3^r))\cong\Z/3^r\oplus\Z/3\oplus\Z/3.\]
  Moreover, it has an order 3 generator $\widetilde{\alpha}_r$ satisfying $q_3^r\circ\widetilde{\alpha}_r\simeq \alpha$.
 
  \begin{proof}
By the homotopy equivalence (\ref{CMN:OmegaF}), we compute that  
\begin{align*}
  \pi_5(\Omega F^3\{3^r\})&\cong \pi_5(S^5\{3^{r+1}\})\oplus \pi_5(\Omega P^4(3^r))\cong \pi_5(\Omega P^4(3^r))\cong \Z/3^r\oplus \Z/3,\\
  \pi_4(\Omega F^3\{3^r\})&\cong \pi_4(S^5\{3^{r+1}\})\oplus \pi_4(\Omega P^4(3^r))\cong \pi_4(S^5\{3^{r+1}\})\cong\Z/3^{r+1},
\end{align*}
where $\pi_5(\Omega P^4(3^r))\cong \Z/3^r\oplus \Z/3$ refers to \cite[Lemma 3.1]{CS}. 
Consider the induced exact sequence 
\begin{multline*}
  0\to \pi_5(\Omega F^3\{3^r\})\to\pi_5(\Omega P^3(3^r))\xra{(\Omega q_3)_\ast}\pi_5(\Omega S^3)_{(3)}\\
    \xra{\partial}\pi_4(\Omega F^3\{3^r\})\xra{(\Omega j_r)_\ast} \pi_4(\Omega P^3(3^r))\to 0.
\end{multline*}
Since $\pi_5(P^3(3^r))\cong\Z/3^{r+1}$ (see \cite[Lemma 2.3]{LZ} or \cite[Theorem 2.10]{Neisen81}), the homomorphism $(\Omega j_r)_\ast$ is an isomorphism, implying that the connecting homomorphism $\partial$ is trivial, and hence there is a short exact sequence
\begin{equation}\label{s.ss}
  0\to \pi_5(\Omega F^3\{3^r\})\xra{(\Omega j_r)_\ast}\pi_5(\Omega P^3(3^r))\xra{(\Omega q_3)_\ast}\pi_5(\Omega S^3)_{(3)}\to 0.
\end{equation}
We claim that this exact sequence splits. For short denote $\Omega\Sigma \bigvee_{\alpha} P^{n_\alpha}(3^r)$ by $\Omega_r$. Since $\pi_5(\Omega F^3\{3^r\})\cong \pi_5(\Omega_r)$, there is a commutative diagram 
\[\begin{tikzcd}
\pi_5(\Omega_r)\ar[r,equal]\ar[d,"(\Psi_3)_\ast","\cong"swap]& \pi_5(\Omega_r)\ar[d,"(\Phi_2)_\ast"]\\
\pi_5(\Omega F^3\{3^r\})\ar[r,"(\Omega j_r)_\ast"]& \pi_5(\Omega P^3(3^r)),
\end{tikzcd}\]
where $\Phi_2$ is the restriction of the homotopy equivalence $\Phi$ in (\ref{cmn-odd}) to $\Omega_r$ and $\Psi_3$ is the restriction of the homotopy equivalence $\Psi$ in (\ref{CMN:OmegaF}) to $\Omega_r$. By \cite[Corollary 11.10.4]{N2}, there is a retraction $\kappa\colon \Omega P^3(3^r)\to \Omega_r$ such that the composition $\kappa\circ(\Phi_2)_\ast$ is homotopic to the identity map.  Therefore $(\Omega j_r)_\ast$ has a retraction and the short exact sequence (\ref{s.ss}) splits. Therefore there is a group isomorphism  
\[\pi_6(P^3(3^r))\cong \Z/3^r\oplus\Z/3\oplus\Z/3.\]
The splitting implies that $(q_3)_\ast\colon \pi_6(P^3(3^r))\to \pi_6(S^3)_{(3)}$ is surjective. Let $\widetilde{\alpha}_r\in\pi_6(P^3(3^r))$ be the generator satisfying $q_3\widetilde{\alpha}_r\simeq\alpha$. This completes the proof. 
  \end{proof}
\end{lemma}

\begin{remark}\label{rmk:pi5T3}
 After localized at $3$, the homotopy fibration sequence (\ref{fib:T}) yields the exact sequence of groups
  \[0\to \pi_5(T^3\{3^r\})\to \pi_5(\Omega S^3)\xra{}\pi_4(S^5\{3^{r+1}\})\to \pi_4(T^3\{3^r\})\to 0.\]
  The proof of Lemma 2.6 implies that $\pi_5(T^3(3^r))\cong \pi_6(S^3)_{(3)}\cong\Z/3$, then the exact sequence above implies that the second arrow, and hence the fourth arrow are group isomorphisms. Thus we have group isomorphisms
\[\pi_5(T^3\{3^r\})\cong \Z/3,\quad \pi_4(T^3\{3^r\})\cong\Z/3^{r+1}.\]

\end{remark}

\begin{lemma}\label{lem:pi7P4}
 There are group isomorphisms 
 \[\pi_7(P^4(3^r))\cong \Z/3\lra{\Sigma\widetilde{\alpha}_r} \text{~ and~ }\pi_7(P^4(p^r))=0  \text{ for $p\geq 5$}.\] 
Moreover, when $p=3$, we have the formula 
\begin{equation}\label{eq:alpha_r}
  B_3(\chi^r_s)\Sigma \widetilde{\alpha}_r\simeq\Sigma\widetilde{\alpha}_s  \text{ for $s\geq r$}
  \end{equation}
\begin{proof}
  For an odd prime $p$, by \cite{CMN79,Neisen81} there is a $p$-local homotopy equivalence \[\Omega P^4(p^r)\simeq_{(p)} S^3\{p^r\} \times \Omega (\vees_{k=0}^{\infty} P^{2k+7}(p^r)). \] 
Hence we have group isomorphisms 
 \[\pi_7(P^4(p^r))\cong \pi_6(S^3\{p^r\})\oplus \pi_7(P^{7}(p^r))\cong\pi_6(S^3\{p^r\})\cong\Z/(3,p).\]
When $p=3$, it remains to show that $\Sigma\widetilde{\alpha}_r$ is a generator of $\pi_7(P^4(3^r))\cong\Z/3$ and the equation (\ref{eq:alpha_r}). By Lemma \ref{lem:pi6p3}, $\widetilde{\alpha}_r$ is an order 3 generator of $\pi_6(P^3(3^r))$ satisfying $q_3\circ \widetilde{\alpha}_r\simeq \alpha$. After suspending we have $q_4\circ \Sigma\widetilde{\alpha}_r\simeq \Sigma\alpha$. Since $\Sigma \alpha$ is nontrivial, so is $\Sigma \widetilde{\alpha}_r$. Hence it has order $3$ and generates $\pi_7(P^4(3^r))$.
 The homotopy (\ref{eq:alpha_r}) follows by (\ref{eq:chi}).
\end{proof}
\end{lemma}

Recall the homotopy cofibration sequence for $C^{n+2}_\eta=\Sigma^{n-2}\mathbb{C}P^2$:
\[S^{n+1}\xra{\eta_n}S^n\xra{i_n^\eta}C^{n+2}_\eta\xra{q_{n+2}^\eta}S^{n+2}.\]

\begin{lemma}\label{lem:piCeta}
 Let $n\geq 3$ and let $1_n$ be the identity map on $S^n$. There hold 
 \begin{align*}
  &\pi_{n}(C^{n+2}_\eta)\cong\Z\lra{i_n^\eta},~~\pi_{n+1}(C^{n+2}_\eta)=0,~~\pi_{n+2}(C^{n+2}_\eta)\cong\Z\lra{\widetilde{\zeta}},\\
  &[C^{n+2}_\eta,S^{n+2}]\cong\Z\lra{q_{n+2}},~~[C^{n+2}_\eta,S^{n+1}]=0,~~[C^{n+2}_\eta,S^n]\cong \Z\lra{\overline{\zeta}},
 \end{align*}
  where $\widetilde{\zeta}$ and $\overline{\zeta}$ satisfy the formulae
 \[q_{n+2}^\eta\widetilde{\zeta}\simeq 2\cdot 1_{n+2},\quad \overline{\zeta}i_n^\eta\simeq 2\cdot 1_n. \]
\begin{proof}
  See \cite[Section 8, Page 92]{Toda2}.
\end{proof}
\end{lemma}

\begin{corollary}\label{cor:hCeta}
  The induced homomorphisms
  \[\widetilde{\zeta}_\ast\colon H_{n+2}(S^{n+2})\to H_{n+2}(C^{n+2}_\eta), \quad \overline{\zeta}_\ast\colon H_n(C^{n+2}_\eta)\to H_n(S^n) \]
  are the multiplication by $2$.
  \begin{proof}
    We only show that $\widetilde{\zeta}_\ast=2$ and omit the proof of $\overline{\zeta}_\ast=2$ which is similar.  Let $\iota_{n+2}\in H_{n+2}(S^{n+2})$ be a generator and let $b_\eta\in H_{n+2}(C^{n+2}_\eta)$ be the class such that $(q_{n+2}^\eta)_\ast(b_\eta)=\iota_{n+2}$.  Set $\widetilde{\zeta}_\ast(\iota_{n+2})=x\cdot b_\eta$ for some integer $x$. Then $q_{n+2}^\eta \widetilde{\zeta}\simeq 2\cdot 1_{n+2}$ implying that 
    \[(q_{n+2}^\eta)_\ast\widetilde{\zeta}_\ast(\iota_{n+2})=(q_{n+2}^\eta)_\ast(x\cdot b_\eta)=x\cdot \iota_{n+2},\] 
     hence $x=2$. 
  \end{proof}
\end{corollary}

\begin{lemma}\label{lem:pi7Ceta}
  There hold isomorphisms
\[\pi_7(C^6_\eta)\cong \Z\lra{i_4^\eta\nu_4}\oplus\Z/6\lra{i_4^\eta(\Sigma\nu')},~\pi_8(C^7_\eta)\cong \Z/12\lra{i_5^\eta\nu_5}\text{ and }\pi_7(C^5_\eta)\cong \Z\lra{\xi},\] 
 where $\xi$ satisfies the formula 
\(\Sigma \xi\simeq i_4^\eta\nu_4\eta_7.\)
\begin{proof}
 See \cite[Propositions 8,2, 8.3 and Lemma 8.5]{Mukai82}.
\end{proof}
\end{lemma}

\subsection{Criteria for null homotopy}
  
The following lemmas will be frequently used in the subsequent sections.

\begin{lemma}\label{lem:Sq2-Ceta}
 For each integer $n\geq 3$, the reduced Steenrod square 
  \[\Sq^2\colon H^n(C^{n+2}_\eta;\z{})\to H^{n+2}(C^{n+2}_\eta;\z{})\] 
  is an isomoprhism. 
  \begin{proof}
   Clear; see also \cite{ZP17}. 
  \end{proof}
\end{lemma}

\begin{lemma}\label{lem:PP}
  Let $C_x$ be the homotopy cofibre of the map $x\alpha\in\pi_6(S^3)$, $xi_3\alpha\in\pi_6(P^4(3^r))$ or $x\widetilde{\alpha}_r\in\pi_6(P^3(3^r))$, where $x\in \{0,\pm 1\}$ and $r\geq 1$. Then the reduced $3$rd power operation
   \[\PP\colon H^3(C_x;\Z/3)\to H^7(C_x;\Z/3)\]
   is an isomorphism if and only if $x=\pm 1$, and is trivial if and only if $x=0$.
   \begin{proof}
There are homotopy commutative diagram 
\[\begin{tikzcd}
  S^6\ar[d,equal]\ar[r,"\widetilde{\alpha}_r"]&P^3(3^r)\ar[r]\ar[d,"q_3"]&C_{\widetilde{\alpha}_r}\ar[d,"\mu"]\\
  S^6\ar[r,"\alpha"]\ar[d,equal]&S^3\ar[d,"i_3"]\ar[r]&C_{\alpha}\ar[d,"\lambda"]\\
  S^6\ar[r,"i_3\alpha"]&P^4(3^r)\ar[r]&C_{i_3\alpha}
\end{tikzcd},\]
where rows are homotopy cofibration sequences and $\mu,\lambda$ are induced maps. It is easy to see that $\mu$ and $\lambda$ induce isomorphisms of mod $3$ cohomology groups in dimensions $3$ and $7$. The Lemma then follows by the naturality of $\PP$ and the fact that $\alpha$ is detected by $\PP$.  
   \end{proof}
   \end{lemma}

\begin{lemma}\label{lem:CS3.3}
  Let $C_y$ be the homotopy cofibre of  $y\cdot [1_P,1_P]\circ \widetilde{i_7}\in \pi_7(P^5(p^r))$, where $y\in\zp{r}$ with $p\geq 3$ an odd prime. Then $y=0$ if and only if all cup products in $H^\ast(C_y;\zp{r})$ are trivial. Moreover, for any $y\in\zp{r}$, the action of the reduced power operation $\PP$ on $H^\ast(C_y;\Z/3)$ is trivial.  
  \begin{proof}
    The proof is similar to that in \cite[Lemma 3.3]{CS}.
  \end{proof}
\end{lemma}

\begin{lemma}\label{lem:CS2.4}
 Let $r,s,m,n,l$ be integers such that  $$
\left\{\begin{array}{c}
r, s \geq 1 \quad \text { and } \quad n, l \geq 3 \\
m=n+l \quad \text { or } \quad m=n+l-1
\end{array}\right.
$$
Let $f\colon S^m \rightarrow \Sigma P^n\left(p^r\right) \wedge P^l\left(p^s\right)$ be a map and let $C_{\hat{f}}$ be the homotopy cofibre of
$$
\hat{f}\colon S^m \xrightarrow{f} \Sigma P^n\left(p^r\right) \wedge P^l\left(p^s\right) \xrightarrow{\left[\imath_1, \imath_2\right]} P^{n+1}\left(p^r\right) \vee P^{l+1}\left(p^s\right),
$$
where $\left[\imath_1, \imath_2\right]$ is the Whitehead product of the inclusions $\imath_1\colon P^{n+1}\left(p^r\right) \rightarrow P^{n+1}\left(p^r\right) \vee P^{l+1}\left(p^s\right)$ and $\imath_2 \colon P^{l+1}\left(p^s\right) \rightarrow P^{n+1}\left(p^r\right) \vee P^{l+1}\left(p^s\right)$. Then the following statements are equivalent:
\begin{enumerate}[(1)]
  \item $f$ is null homotopic;
  \item$f^\ast\colon H^\ast\left(\Sigma P^n\left(p^r\right) \wedge P^l\left(p^s\right) ; \mathbb{Z} / p^{m_{r,s}}\right) \rightarrow H^*\left(S^m ; \mathbb{Z} / p^{m_{r,s}}\right)$ is trivial, where $m_{r,s}=\min\{r,s\}$;
 \item  all cup products in $\tilde{H}^*\left(C_{\hat{f}} ; \mathbb{Z} / p^t\right)$ are trivial.
\end{enumerate}
Moreover, the Lemma holds if we allow exactly one of $r,s$ to be $\infty$ under the convention $P^{n+1}(p^{\infty})=S^n$. 
\begin{proof}
  See \cite[Lemma 2.4]{CS}.
\end{proof}
\end{lemma}

Let $p$ be a prime and let $\beta_r\colon H^\ast(-;\Z/p)\dashrightarrow H^{\ast+1}(-;\Z/p)$ be the $r$-th higher order Bockstein operation defined in the Introduction. 
\begin{lemma}\label{lem:Bockstein}
  The following statements hold:
    \begin{enumerate} 
      \item The higher Bockstein $\beta_r$ detects the degree $p^r$ map on $S^n$; in other words,  for each $r\geq 1$, there are exactly one non-trivial higher Bockstein  
      \[\beta_r\colon H^{n-1}(P^n(p^r);\zp{})\to H^n(P^n(p^r);\zp{}).\]
      \item For each $r\geq 1$, elements of $H^\ast(X;\zp{})$ coming from free integral homology class lie in $\ker(\beta_r)$ and not in $\im(\beta_r)$.
    \item If $z\in H^{n+1}(X;\Z)$ generates a direct summand $\zp{r}$ for some $r$, then there exist generators  $z'\in H^n(X;\zp{})$ and $z''\in H^{n+1}(X;\zp{})$ such that 
      \[ \beta_r(z')=z'',\quad \beta_i(z')=\beta_i(z'')=0 \text{~for~} i<r.\]
    \end{enumerate}
    \begin{proof}
      The first statement is well-known, compare \cite[Page 173, Proposition 1]{MT68} for $p=2$; the second and third statements for $p=2$ refer to \cite[Page 61, Proposition 2]{MT68}, and the proofs of the statements for odd primes $p$ are similar.
    \end{proof}
  \end{lemma}

\begin{lemma}[compare Lemma 5.6 of \cite{ST}]\label{STlemma}
 Consider a homotopy cofibration sequence of simply connected CW-complexes:
\[
S\xra{f} \vees_{i} A_i\vee B\stackrel{g}{\longrightarrow} \Sigma C.
\]
If the composition $S\xra{f} \vees_{i} A_i\vee B\xra{q_j}A_j$, where $q_j$ is the canonical projection, is null homotopic for each $j$, then there is a homotopy equivalence
\[
\Sigma C\simeq \vees_{i} A_i\vee D,
\]
where $D$ is the homotopy cofibre of the composition $S\xra{f} \vees_{i} A_i\vee B
\xra{q_B} B$ with $q_B$ the projection .
\end{lemma} 
\begin{proof}
Consider the diagram of homotopy cofibration sequences
\[\begin{tikzcd}
  S\ar[r,"f"]\ar[d,equals]&\bigvee_{i} A_i\vee B\ar[r,"g"]\ar[d,"q_i"]&\Sigma C\ar[d,"h_i"]\\
  S\ar[r,"f_i"]&A_i\ar[r,"g_i"]&F_i,
\end{tikzcd}\]
where $q_i$ is the projection onto $A_i$, $F_i$ is the homotopy cofibre of $f_i= q_i\circ f$, and $g_i,h_i$ are induced maps. By the condition $f_i$ is null homotopic, and thus $F_i\simeq A_i\vee \Sigma S$. In particular, $g_i$ has a left homotopy inverse $\varrho_i$.
It means that the restriction of $g$ on each $A_i$ has a left homotopy inverse $\varrho_i\circ h_i$.
Define the composition 
\[
\varrho\colon \Sigma C \stackrel{\mu^\prime}{\longrightarrow} \vees_{i}\Sigma C \xra{\vees_{i}h_i}\vees_{i}F_i \xra{\vees_{i}\varrho_i} \vees_{i}A_i,
\]
where $\mu^\prime$ is the iterated comultiplication of $\Sigma C$. Then the homotopy commutative diagram 
\[\begin{tikzcd}
\Sigma C\ar[dr,equal]\ar[r,"\mu'"]&\vees_i\Sigma C\ar[d,"p_i"]\ar[r,"\vees_i h_i"]&\vees_i F_i\ar[r,"\vees_i\varrho_i"]\ar[d,"p_i"]&\vees_i A_i\ar[d,"p_i"]\\
&\Sigma C\ar[r,"h_i"]&F_i\ar[r,"\varrho_i"]&A_i
\end{tikzcd}\]
shows that there is a homotopy $p_i\circ \varrho\simeq \varrho_i h_i$ for each $i$, where $p_i$ is the pinch map collasping the subspace $\vees_{j\neq i}A_j$ to the base point,  $1\leq i,j\leq n$.
We \emph{claim} that $\varrho$ is a left homotopy inverse of $g_A\colon A=\vees_{i}A_i \stackrel{i_A}{\longrightarrow}   \vees_{i} A_i\vee B \stackrel{g}{\longrightarrow}   \Sigma C$, where $i_A$ is the inclusion map. Denote by $\iota_j$ the canonical inclusion $A_j\to\vees_i A_i$, $j=1,\cdots,n$.	Consider the following homotopy commutative diagram 
\[\begin{tikzcd}
\vees_i A_i\ar[r,"i_A"]&\vees_{i} A_i\vee B\ar[d,"q_i"]\ar[r,"g"]&\Sigma C\ar[d,"h_i"]\ar[r,"\varrho"]&\vees_iA_i\ar[d,"p_i"]\\
A_j\ar[u,"\iota_j"]\ar[r,"\delta_{ij}\cdot id_{A_j}"]&A_i\ar[r,"g_i"]&F_i\ar[r,"\varrho_i"]&A_i 
\end{tikzcd},\]
where $\delta_{ij}=1$ for $i=j$, and $\delta_{ij}=0$ for $i\neq j$. It follows that 
\begin{align*}
  p_i\circ\varrho\circ g_A\circ \iota_j\simeq \varrho_i\circ g_i\circ (\delta_{ij}\cdot id_{A_j})\simeq \delta_{ij}\cdot id_{A_j}, ~\forall~1\leq i,j\leq n.
\end{align*}
Passing to the induced homomorphism of integral homology, we see that the homomorphism of integral homology groups induced by the composition $\varrho\circ g_A$ is the identity, or the diagonal matrix with identity	diagonal entries.  Then by the Whitehead's Theorem, we get that $\varrho\circ g_A$ is a self-homotopy equivalence of $A$. Replacing $\varrho$ by $(\varrho\circ g_A)^{-1}\varrho$ if necessary, we get $\varrho \circ g_A\simeq id_A$, and thus the claim is proved.

Consider the homotopy commutative diagram of homotopy cofibrations
\[\begin{tikzcd}
  \ast\ar[d]\ar[r]&\vees_{i} A_i\ar[r,equals]\ar[d,"i_A"]&\vees_{i}A_i\ar[d,"g_A"]\\
  S\ar[r,"f"]\ar[d,equals]&\vees_{i} A_i\vee B\ar[r,"g"]\ar[d,"q_B"]&\Sigma C\ar[d,"q"]\\
  S\ar[r,"q_B\circ f"]&B\ar[r,"q'"]&D
\end{tikzcd},\]
where $q$ and $q'$ are induced maps. The third column is a homotopy cofibration sequence in which $g_A$ admits a left homotopy inverse $\varrho$, so we have a homotopy equivalence 
\[\varrho+q\colon \Sigma C\xra{\mathrm{comulti}}\Sigma C\vee \Sigma C\xra{\varrho\vee q}\vees_iA_i\vee D.\]
The proof of the Lemma is finished.
\end{proof}

\section{Homology sections of $\Sigma M$}
\label{sec: hdec}
Let $M$ be a closed simply connected $7$-manifold with $H_\ast(M;\Z)$ given by (\ref{HMeq}). 
  By the homology groups $H_\ast(M;\Z)$ (\ref{HMeq}) and the homology decomposition theorem (see \cite[Theorem 4H.3]{Hat}), there are homotopy cofibration sequences:
\begin{equation}\label{Cof:M}
  \begin{aligned}
    &\bigvee_{i=1}^kS^2 \vee P^{3}(T)\stackrel{f_{2}}{\longrightarrow} M_2\to M_3, ~ \ {\rm where}~ 
     M_2\simeq\bigvee_{i=1}^lS^2 \vee P^3(H);\\
    & \bigvee_{i=1}^kS^3 \vee P^4(H)\stackrel{f_{3}}{\longrightarrow} M_3\to M_4,~~\bigvee_{i=1}^lS^4\stackrel{f_{4}}{\longrightarrow}  M_4\to M_5,\\
    &S^6\xra{f_5}M_5\to M,
  \end{aligned} 
\end{equation}
where the attaching maps $f_2,f_3,f_4,f_5$ are \emph{homologically trivial} (which means they induce trivial homomorphisms in homology). The complexes $M_i$ are called the \emph{$i$-th homology section} of $M$ and we have 
\[H_j(M_i)\cong H_j(M) \text{ for $j\leq i$ and }H_j(M_i)=0 \text{ for $j>i$}.\]
It is clear that the suspension of the above homotopy cofibration sequences give the homology decomposition of $\Sigma M$. In the remainder of this section we assume that \emph{the torsion groups $H$ and $T$ are both $2$-torsion-free}. 



\begin{lemma}\label{lem:SM3}
There is a homotopy equivalence 
  \begin{equation*}
    \Sigma M_3\simeq \bigvee_{i=1}^{l}S^3\vee \bigvee_{j=1}^{k}S^4 \vee P^4(H) \vee P^5(T).
  \end{equation*}
  \begin{proof}
By (\ref{Cof:M}),  there is a homotopy cofibration sequence 
  \[\bigvee_{i=1}^{k}S^3 \vee P^{4}(T)\xra{\Sigma f_2}\Sigma M_2\to\Sigma M_3, \]
  where $\Sigma f_2$ is homologically trivial and $\Sigma M_2\simeq \bigvee_{i=1}^{l}S^3 \vee P^4(H)$. 
Consider the following homologically trivial compositions
\begin{align*}
&S^3\hookrightarrow\bigvee_{i=1}^{k}S^3 \vee P^{4}(T)\xra{\Sigma f_2}\Sigma M_2\twoheadrightarrow X \text{ for }X\in\{S^3,P^4(p^s)\};\\
&P^4(p^r)\hookrightarrow P^{4}(T)\hookrightarrow\bigvee_{i=1}^{k}S^3 \vee P^{4}(T)\xra{\Sigma f_2}\Sigma M_2\twoheadrightarrow Y\text{ for }Y\in\{S^3,P^4(q^s)\},
\end{align*}  
  where the unlabeled arrows ``$\hookrightarrow$'' and ``$\twoheadrightarrow$'' denote the canonical inclusion and pinch maps, respectively. Note that Lemma \ref{lem:Moore1} (\ref{PrPs}) implies that for each $n\geq 4$ and each odd prime $p$, there is a bijection 
  \begin{equation}\label{eq:PrPs}
    [P^n(p^r),P^n(p^s)]\xra{}\Hom\big(H_{n-1}(P^n(p^r)),H_{n-1}(P^n(p^s))\big),\quad B_{n-1}(\chi^r_s)\mapsto \chi^r_s,
  \end{equation} 
  where $\chi^r_s$ refers to (\ref{eq:chi^r_s}); 
  see also \cite[Proposition 1.4.2]{N2}. In particular, a map $P^n(p^r)\to P^n(p^s)$ is null homotopic if and only if it is homologically trivial. Hence by the Hurewicz theorem and Lemma \ref{lem:Moore1}, we see that all the above compositions are null homotopic. The Lemma then follows by Lemma \ref{STlemma}.
  \end{proof}
\end{lemma}

\begin{lemma}\label{lem:SM4}
There is a homotopy equivalence 
  \begin{equation*}
    \Sigma M_4\simeq P^4(H) \vee P^5(T)  \vee P^6(H)\vee \bigvee_{i=1}^{l-a}S^3\vee \bigvee_{j=1}^{k}S^4 \vee  \bigvee_{r=1}^{k-a}S^5\vee  \bigvee_{s=1}^aC^5_\eta,
  \end{equation*} 
  where $0\leq a\leq \min\{k,l\}$ is some integer that depends on $M$. Moreover, $a=0$ if and only if the Steenrod square $\Sq^2$ acts trivially on $H^2(M;\z{})$.
  \begin{proof}
   By (\ref{Cof:M}) and  Lemma \ref{lem:SM3},  there is a homotopy cofibration sequence 
  \[ \bigvee_{i=1}^kS^4 \vee P^5(H)\xra{\Sigma f_3} \Sigma M_3\to \Sigma M_4,\]
  where $\Sigma f_3$ is homologically trivial and $\Sigma M_3\simeq \bigvee_{i=1}^{l}S^3\vee \bigvee_{j=1}^{k}S^4 \vee P^4(H) \vee P^5(T)$.
  By the Hurewicz theorem, Lemma \ref{lem:Moore1} and (\ref{eq:PrPs}),  the following compositions
  \begin{align*}
    &S^4\hookrightarrow \bigvee_{j=1}^kS^4 \vee P^5(H)\xra{\Sigma f_3} \Sigma M_3\twoheadrightarrow X\text{ for }X\in \{S^4, P^4(p^s),P^5(p^t)\};\\
    &P^5(p^s)\hookrightarrow P^5(H)\hookrightarrow \bigvee_{j=1}^kS^4 \vee P^5(H)\xra{\Sigma  f_3} \Sigma M_3\twoheadrightarrow Y \text{ for } Y\in\{S^3,S^4,P^4(q^r), P^5(q^t)\}
  \end{align*}
  are null homotopic, where the unlabelled maps are the canonical inclusions and pinch maps, respectively. 
  It follows by Lemma \ref{STlemma} that there is a homotopy equivalence 
  \[\Sigma M_4\simeq \bigvee_{j=1}^{k}S^4 \vee P^4(H) \vee P^5(T)\vee P^6(H)\vee C_{g_1},\]
  where $C_{g_1}$ is the homotopy cofibre of some map $g_1\colon \bigvee_{j=1}^{k}S^4 \to \bigvee_{i=1}^{l}S^3$.
  The map $g_1$ can be represented by the mod-2 matrix $M_{g_1}=[m_{ij}\eta]_{k\times l}$, where $m_{ij}\in\Z/2$ and the $(i,j)$-entry $m_{ij}\eta$ is the composition $S^4\hookrightarrow \bigvee_{j=1}^kS^4\xra{g_1}\bigvee_{i=1}^lS^3\twoheadrightarrow S^3$ of the inclusion of the $i$-th $4$-sphere into the wedge and the projection onto the $j$-th $3$-sphere. Note that applying an elementary row operation to $M_{g_1}$ corresponds to precomposing a homotopy self-equivalence of $\bigvee_{j=1}^k S^4$ to $g_1$, and the mapping cone of the resulting map is homotopy equivalent to $C_{g_1}$. Similarly, the homotopy type of $C_{g_1}$ remains unchanged if we apply an elementary column operation to $M_{g_1}$. Using a sequence of suitable elementary matrix operations, we may assume that $M_{g_1}$ is a (rectangular) diagonal matrix. Let $a$ be the rank of $M_{g_1}$. Then there is a homotopy equivalence
  \[C_{g_1}\simeq \bigvee_{i=1}^{l-a}S^3\vee \bigvee_{j=1}^{k-a}S^5\vee \bigvee_{s=1}^aC^5_\eta.\]
 Thus we get the homotopy equivalence for $\Sigma M_4$ in the Lemma. 
 
 By (\ref{Cof:M}) and the universal coefficient theorem for cohomology, we have 
 \[H^i(M;\z{})\cong H^i(M_4;\z{}) \text{ for }i\leq 4.\] 
It follows that $\Sq^2$ acts trivially on $H^2(M;\z{})$ if and only if $\Sq^2$ acts trivially on $H^2(M_4;\z{})$. By the homotopy equivalence for $\Sigma M_4$ in the Lemma, we see that $\Sq^2$ acts trivially on $H^2(M;\z{})$ if and only if $\Sq^2$ acts trivially on $H^3(\bigvee_{s=1}^aC^5_\eta;\z{})$, which happens if and only if $a=0$ by Lemma \ref{lem:Sq2-Ceta}.
  The proof of the Lemma is finished.
  \end{proof}
  
\end{lemma}

Note that if $M$ is a smooth spin manifold, then similar arguments to that on \cite[page 32]{MM79} show that  the secondary operation $\Theta$ that detects $\eta^2\in\pi_{n+2}(S^n)$ acts trivially on $H^\ast(M;\z{})$; see also \cite[Lemma 4.3]{LZ}.

\begin{lemma}\label{lem:SM5-glob}
If $\Theta$ acts trivially on $H_\ast(M;\z{})$, then there is a homotopy equivalence 
  \[
    \Sigma M_5\simeq  P^4(H)\vee P^5(T)\vee P^6(H)\vee \bigvee_{i=1}^{l-a}S^3\vee \bigvee_{j=1}^{k-b}S^4\vee \bigvee_{r=1}^{k-a}S^5\vee \bigvee_{s=1}^{l-b}S^6\vee \bigvee_{u=1}^a C^5_\eta\vee \bigvee_{v=1}^bC^6_\eta,
    \]
    where $0\leq a,b\leq \min\{k,l\}$. Moreover, $a=0$ if and only if $\Sq^2$ acts trivially on $H^2(M;\z{})$, and $b=0$ if and only if $\Sq^2$ acts trivially on $H^3(M;\z{})$.
    \begin{proof}
 By (\ref{Cof:M}), there is a homotopy cofibration sequence
     \[\bigvee_{i=1}^{l}S^5\xra{\Sigma f_4}\Sigma M_4\to\Sigma M_5,\]
     where $\Sigma f_4$ is homologically trivial and $\Sigma M_4$ is given by Lemma \ref{lem:SM4}. Note that the compositions 
     \[S^5\hookrightarrow \bigvee_{i=1}^{l}S^5\xra{\Sigma f_4}\Sigma M_4\hookrightarrow Y\]
     are null homotopic for $Y=S^5,P^4(H),P^5(T),P^6(T)$.
     Since the induced homomorphism $\widetilde{\zeta}_\ast\colon H_5(S^5)\to H_5(C^5_\eta)$ is the multiplication by $2$ (Corollary \ref{cor:hCeta}), we see that the homologically trivial composition
     \[S^5\hookrightarrow \bigvee_{i=1}^{l}S^5\xra{\Sigma f_4}\Sigma M_4\hookrightarrow \bigvee_{i=1}^a C^5_\eta\hookrightarrow C^5_\eta\]
     is null homotopic. 
     By  Lemmas \ref{lem:Moore1} and \ref{STlemma}, we see that there is a homotopy equivalence
     \begin{equation}\label{SM5-glob}
      \Sigma M_5\simeq \bigvee_{i=1}^{k-a}S^5\vee P^4(H)\vee P^5(T)\vee P^6(H)\vee \bigvee_{j=1}^a C^5_\eta\vee C_{g_2},
     \end{equation} 
     where $C_{g_2}$ is the homotopy cofibre of some map 
     \[g_2\colon \bigvee_{i=1}^{l}S^5\to \bigvee_{r=1}^{l-a}S^3\vee \bigvee_{s=1}^kS^4.\]
     By assumption, the secondary operation $\Theta$ that detects $\eta^2$ acts trivially on $H^\ast(M;\z{})$, and hence acts trivially on $H^\ast(C_{g_2};\z{})$. It follows that the map $g_2$ contains no $\eta^2$ components, and hence the space $ \bigvee_{i=1}^{l-a}S^3$ retracts off $C_{g_2}$. By similar arguments to that in the proof of Lemma \ref{lem:SM4}, we see that there is a homotopy equivalence
     \[C_{g_2}\simeq  \bigvee_{i=1}^{l-a}S^3\vee \bigvee_{j=1}^{k-b}S^4\vee \bigvee_{r=1}^{l-b}S^6\vee \bigvee_{s=1}^bC^6_\eta\] 
     for some integer $0\leq b\leq \min\{k,l\}$. Thus we obtain the homotopy equivalence for $\Sigma M_5$ in the Lemma.

    The proof of the second part is similar to that given at the last paragraph of the proof of Lemma \ref{lem:SM4}.
    \end{proof}
\end{lemma}

\begin{corollary}\label{cor:SM5-glob}
   If M is spin and smooth, then there is a homotopy equivalence 
   \[
    \Sigma M_5\simeq  P^4(H)\vee P^5(T)\vee P^6(H)\vee \bigvee_{i=1}^{l-a}S^3\vee \bigvee_{j=1}^{k-b}S^4\vee \bigvee_{r=1}^{k-a}S^5\vee \bigvee_{s=1}^{l-b}S^6\vee \bigvee_{u=1}^a C^5_\eta\vee \bigvee_{v=1}^bC^6_\eta,
    \]
    where $0\leq a,b\leq \min\{k,l\}$. Moreover, $a=0$ if and only if $\Sq^2$ acts trivially on $H^2(M;\z{})$, and $b=0$ if and only if $\Sq^2$ acts trivially on $H^3(M;\z{})$.
\end{corollary}

\section{Homotopy decompositions of $\Sigma M$ when localized away from $2$}
\label{sec:SM:p-local}
 Let $M$ be a closed simply connected $7$-manifold with $H_\ast(M;\Z)$ given by (\ref{HMeq}), where $H$ and $T$ are both $2$-torsion-free finite abelian groups. In this section we discuss the homotopy type of the reduced suspension $\Sigma M$ after localization away from $2$.
Firstly we have the following immediate consequence of Lemma \ref{lem:SM5-glob}.

\begin{lemma}\label{lem:SM5}
There is a homotopy equivalence 
\[\Sigma M_5\simeq_{\loca} \bigvee_{i=1}^{l}(S^3\vee S^6)\vee \bigvee_{j=1}^{k}(S^4\vee S^5)\vee P^4(H)\vee P^5(T)\vee P^6(H).\]
\end{lemma}

Set $H=H_3\oplus H_{\neq 3}$ and $T=T_3\oplus T_{\neq 3}$ with \[H_3=\bigoplus_{i=1}^{m}\Z/3^{r_i},\quad T_3=\bigoplus_{j=1}^n\Z/3^{s_j},\]
where $H_{\neq 3}$ and $T_{\neq 3}$ are subgroups of $H$ and $T$ with no $3$-torsion. 
\begin{lemma}\label{lem:SM}
  There is a homotopy equivalence
  \[
  \Sigma M\simeq_{\loca} \bigvee_{i=1}^l(S^3\vee S^6)\vee \bigvee_{j=1}^k S^5  \vee P^6(H) \vee  P^4(H_{\neq 3}) \vee P^5(T_{\neq 3})\vee C_{h},
  \]
  where $C_h$ is the homotopy cofibre of the map \[h\colon S^7\xra{\Sigma f_5}\Sigma M_5\twoheadrightarrow \bigvee_{j=1}^{k}S^4\vee P^4(H_3) \vee P^5(T_3).\]
  \begin{proof}
    Consider the cofibration $S^7\xra{\Sigma f_5}\Sigma M_5\to \Sigma M$. We \emph{claim} that when localized away from $2$, the compositions 
    \[S^7\xra{\Sigma f_5}\Sigma M_5\twoheadrightarrow X\]
    are null homotopic for $X=S^3,S^6, S^5, P^4(H_{\neq 3}),P^5(T_{\neq 3}),P^6(H)$, where ``$\twoheadrightarrow$'' denotes the canonical pinch map by Lemma \ref{lem:SM5}. If so, the Lemma follows by Lemma \ref{STlemma}. 

    For $X=S^3,S^6, S^5,P^6(H)$, it is due to the fact that $\pi_7(S^3)$, $\pi_7(S^6)$, $\pi_7(S^5)$ and $\pi_7(P^6(H))$ are trivial when  localized away from $2$, by \cite{Tod} and Lemma \ref{lem:Moore1}. For $X=P^4(H_{\neq 3})$,  we have $\pi_7(P^4(p^r))=0$ for a prime $p\geq 5$ by Lemma \ref{lem:pi7P4}. For $X=P^5(T_{\neq 3})$, Lemma \ref{lem:pi7P5} implies that for primes $p\geq 5$, the composition of the form $$S^7\xra{\Sigma f_5}\Sigma M_5 \twoheadrightarrow P^5(p^r)$$ 
    is given by $f_x=x\cdot [1_P,1_P]\circ \widetilde{i_7}$ for some $x\in\zp{r}$. 
   Consider the homotopy commutative diagram \[\begin{tikzcd}
   S^7\ar[r,"\Sigma f"]\ar[d,equal]&\Sigma M_5\ar[d]\ar[r]&\Sigma M\ar[d,"\lambda"]\\
S^7\ar[r,"f_x"]&P^5(p^r)\ar[r]&C_{f_x}
   \end{tikzcd},\]
   where $\lambda$ is an induced map. Note that the induced homomorphism 
   \[\lambda^\ast\colon H^\ast(C_{f_x};\zp{r})\to H^\ast(\Sigma M;\zp{r})\] 
   is a mononomorphism of cohomology rings.
    Then all cup products in $H^\ast(\Sigma M;\zp{r})$ are trivial implying that all cup products in $H^\ast(C_{f_x};\zp{r})$ are trivial, which in turn implies that $x=0$ by Lemma \ref{lem:CS3.3}. 
  \end{proof}
\end{lemma}

Using the homotopy equivalences $P^4(H_3)\simeq \bigvee_{u=1}^m P^4(3^{r_u})$ and $P^5(T_3)\simeq \bigvee_{v=1}^{n}P^5(3^{s_v})$ and the Hilton-Milnor theorem, we let 
\[V\coloneqq\bigvee_{i=1}^{k}S^4_i\vee \bigvee_{u=1}^m P^4(3^{r_u}) \vee \bigvee_{v=1}^{n}P^5(3^{s_v}),\]
where $S^4_i$ is the $i$-th wedge summand, and express the map $h\colon S^7\to V$ in Lemma \ref{lem:SM} as 
\[h\simeq M_h+\theta_h,\]
where $\theta_h$ is a sum of Whitehead products, $M_h$ is the composite
\[S^7\xra{\mu'_S}\vees_{i,u,v}S^7\xra{\vees_i h_{S,i}\vee \vees_{u}h_{H,u}\vee \vees_vh_{T,v}}\bigvee_{i=1}^{k}S^4_i\vee \bigvee_{u=1}^m P^4(3^{r_u}) \vee \bigvee_{v=1}^{n}P^5(3^{s_v})=V\]
with $\mu'_S$ being is the standard comultiplication on $S^7$, and
\begin{align*}
  h_{S,i}\colon &S^7\xra{h} V\twoheadrightarrow S^4_i,~~i=1,\cdots,k;\\
  h_{H,u}\colon &S^7\xra{h} V\twoheadrightarrow P^4(3^{r_u}),~~u=1,\cdots,m;\\
  h_{T,v}\colon &S^7\xra{h} V\twoheadrightarrow P^5(3^{s_v}),~~v=1,\cdots,n.
\end{align*}
 We abuse the notation slightly and write $M_h$ as the \emph{formal sum}
\[M_h\simeq \sum_{i=1}^k h_{S,i}+ \sum_{u=1}^mh_{H,u}+\sum_{v=1}^nh_{T,v}\]
for convenience and to emphasize the components $h_{S,i},h_{H,u}$ and $h_{T,v}$.  
Alternatively, we may write $M_h$ as the \emph{representation vector} 
\[M_h=[h_{S,1},\ldots,h_{S,k},h_{H,1},\ldots,h_{H,m}, h_{T,1},\ldots,h_{T,n}]^t,\]
where $[-,\cdots,-]^t$ denotes the transpose of the row vector. Here we introduce the notation of representation vector $M_h$ because it is convenient to apply elementary matrix operations, which means composing a self-homotopy equivalence to $h$ and does not change the homotopy type of $C_h$. Consequently, the reduction of the representation vector $M_h$ helps to decompose the homotopy cofibre $C_h$ into simple pieces up to homotopy when $\theta_h=0$.

The cohomology ring $H^\ast(\Sigma M;R)$ has trivial cup products for any principle integral domain $R$ implying that so does all cup products in $H^\ast(C_h;R)$ by Lemma \ref{lem:SM}. Then it follows by \cite[Lemma 4.2]{ST} that the cup products in $H^\ast(C_{h_{S,i}};\Z)$, $H^\ast(C_{h_{H,u}};\Z/3^{r_u})$ and $H^\ast(C_{h_{T,v}};\Z/3^{s_v})$ are all trivial. Since $\pi_7(S^4)_{\loca}\cong \Z_{\loca}\lra{\nu_4}\oplus \Z/3\lra{\Sigma\alpha}$ and the Hopf map $\nu_4$ has Hopf invariant $1$, we apply \cite[Lemma 2.5 (2)]{CS} to get 
\[h_{S,i}\simeq x_i\cdot \Sigma\alpha, \text{  where }x_i\in\Z/3.\]
By Lemma \ref{lem:pi7P4}, $h_{H,u}\simeq y_u\cdot \Sigma\widetilde{\alpha}_{r_u}$, where $y_u\in\Z/3$.
By Lemmas \ref{lem:pi7P5}, \ref{lem:CS3.3} and \cite[Lemma 2.5 (2)]{CS}, we have 
\[h_{T,v}\simeq z_v\cdot i_4(\Sigma \alpha), \text{  where }z_v\in\Z/3.\] 
By the Hilton-Milnor theorem, the Whitehead product component $\theta_{h}$ is determined by groups of the following form:
\begin{align*}
  &\pi_7(\Sigma S^3_i\wedge S^3_j),~~\pi_7(\Sigma S^3_i\wedge P^4(3^{s_j})),~~\pi_7(\Sigma P^3(3^{r_i})\wedge P^4(3^{s_j})),~~\pi_7(\Sigma P^4(3^{s_i})\wedge P^4(3^{s_j})).
\end{align*}
Using the convention $S^n=P^{n+1}(3^{\infty})$, these groups can be written uniformly as 
\[\pi_7(\Sigma P^{a}(3^{r_i})\wedge P^b(3^{s_j})), \text{ where $a+b=7$ or $a+b=8$}.\] 
That is, $\theta_h$ is a sum of maps of the form 
\begin{equation}\label{eq:theta_h}
 \theta_{ij}\colon S^7\to \Sigma P^a(3^{r_i})\wedge P^b(3^{s_j})\xra{[\imath_1,\imath_2]}\Sigma P^a(3^{r_i})\vee \Sigma P^b(3^{s_j})\hookrightarrow V.
\end{equation} 
By \cite[Lemma 2.5]{CS}, all cup products in $H^\ast(\Sigma M;R)$ are trivial implying that all cup products in $H^\ast(C_{\theta_h};R)$ are trivial; this also holds when $\theta_{h}$ is replaced by its component $\theta_{ij}$ of the form (\ref{eq:theta_h}): Let $h_{ij}$ be the composite
\[S^7\xra{h}V\xra{\mathrm{pinch}}\Sigma P^a(3^{r_i})\vee \Sigma P^b(3^{s_j})\hookrightarrow V.\]
Then $h_{ij}\simeq h_i+h_j+\theta_{ij}$, where $h_i$ and $h_j$ are the components of $h_{ij}$ corresponding to $\Sigma P^a(3^{r_i})$ and $\Sigma P^b(3^{s_j})$, respectively. By \cite[Lemma 4.2]{ST}, all cup products in $H^\ast(C_{h_{ij}};R)$ are trivial. As shown above, all cup products in $H^\ast(C_{h_i};R)$ and $H^\ast(C_{h_j};R)$. Hence by \cite[Lemma 2.5 (2)]{CS}, all cup products in $H^\ast(C_{\theta_{ij}};R)$ are also trivial.
 Applying  Lemma \ref{lem:CS2.4}, we see that all such maps are null homotopic. Thus $\theta_h=0$.

Summarizing the above discussion we get 
\begin{lemma}\label{lem:eq-h}
  The map $h\colon S^7\to \bigvee_{i=1}^{k}S^4\vee \bigvee_{u=1}^mP^4(3^{r_u}) \vee \bigvee_{v=1}^{n}P^5(3^{s_v})$ in Lemma \ref{lem:SM} is given by the equation
  \begin{equation}\label{eq:h}
    h\simeq\sum_{i=1}^kx_i\cdot \Sigma\alpha+\sum_{u=1}^my_u\cdot \Sigma\widetilde{\alpha}_{r_u}+\sum_{v=1}^nz_v\cdot i_4(\Sigma \alpha),
  \end{equation}
  where $x_i,y_u,z_v\in \{0,\pm 1\}$. In particular, the homotopy cofibre $C_h$ is a suspension.
\end{lemma}

Note that the map $h$ in (\ref{eq:h}) can be written as the representation vector 
\[M_h=[x_1\cdot \Sigma \alpha,\ldots,x_k\cdot \Sigma\alpha,y_1\cdot \Sigma\widetilde{\alpha}_{r_1},\ldots,y_m\cdot\Sigma\widetilde{\alpha}_{r_m},z_1\cdot i_4(\Sigma \alpha),\ldots,z_n\cdot i_4(\Sigma\alpha)]^t.\]

\begin{proposition}\label{prop:PP=0}
  If the reduced power operation $\PP$ acts trivially on $H^4(C_h;\Z/3)$, then there is a homotopy equivalence 
  \[C_h\simeq  \vees_{i=1}^{k}S^4\vee \bigvee_{u=1}^mP^4(3^{r_u}) \vee \bigvee_{v=1}^{n}P^5(3^{s_v})\vee S^8.\]
\begin{proof}
 By Lemma \ref{lem:PP}, the assumption compels that the coefficients $x_i,y_u,z_v$ in (\ref{eq:h}) are zero for any $i,u,v$. Hence the map $h$ is null homotopic and the Lemma follows. 
\end{proof}
\end{proposition}

\begin{proposition}\label{prop:PP=1}
  Suppose that $\PP$ acts non-trivially on $H^4(C_h;\Z/3)$. The following hold:
  \begin{enumerate}
    \item\label{PP=1:P4} If there exists $x\in H^4(C_h;\Z/3)$ such that $\PP(x)\neq 0$ and $x\in \im (\beta_{r})$ for some $r$, then there is a homotopy equivalence 
    \[C_h\simeq \bigvee_{i=1}^{k}S^4\vee P^4\big(\frac{H_3}{\Z/3^{\bar{r}}}\big)\vee P^5(T_3)\vee (P^4(3^{\bar{r}})\cup_{\Sigma \widetilde{\alpha}_{\bar{r}}}e^8),\]
    where $\bar{r}$ is the minimum of $r$ such that $x\in \im(\beta_r)$ and $\PP(x)\neq 0$.
    \item\label{PP=1:S4} If for any $x\in H^4(C_h;\Z/3)$ with $\PP(x)\neq 0$, there hold $x\notin \im(\beta_r)$ and $\beta_s(x)=0$ for any $r,s\geq 1$, then 
    there is a homotopy equivalence 
    \[C_h\simeq \bigvee_{i=1}^{k-1}S^4\vee P^4(H_3)\vee P^5(T_3)\vee (S^4\cup_{\Sigma\alpha}e^8).\]
    \item\label{PP=1:P5} If for any $x\in H^4(C_h;\Z/3)$ with $\PP(x)\neq 0$, there hold  $x\notin \im(\beta_r)$ for any $r\geq 1$, while $\beta_s(x)\neq 0$ for some $s\geq 1$, then there is a homotopy equivalence 
    \[C_h\simeq \bigvee_{i=1}^{k}S^4\vee P^4(H_3)\vee P^5\big(\frac{T_3}{\Z/3^{\bar{s}}}\big)\vee (P^5(3^{\bar{s}})\cup_{i_4(\Sigma\alpha)}e^8),\]
    where $\bar{s}$ is the maximum of $s$ such that $\beta_s(x)\neq 0$.
  \end{enumerate}
  
  \begin{proof}
  By Lemma \ref{lem:PP},  $\PP\big(H^4(C_h;\Z/3)\big)\neq 0$ implying that there is at least one of $x_i,y_i,z_j$ in (\ref{eq:h}) is nonzero. Since the attaching map $h$ is a suspension, we can use the matrix method in \cite[Section 3]{Hua2} or \cite[Section 2]{LZ} to reduce the coefficients of the equation (\ref{eq:h}).

  Firstly, by Lemma \ref{lem:Bockstein}, the assumption of (\ref{PP=1:P4}) implies that $y_u=\pm 1$ for some $1\leq u\leq m$.  The formulae (\ref{eq:alpha_r}) implies the following homotopies:
  \begin{align*}
    &\mat{1_P}{0}{-q_4^r}{1_P}\matwo{\Sigma\widetilde{\alpha}_r}{\Sigma\alpha}\simeq\matwo{\Sigma\widetilde{\alpha}_r}{0}\colon S^7\to P^4(3^r)\vee S^4,\\
    &\mat{1_P}{0}{-B_3(\chi^r_s)}{1_P}\matwo{\Sigma\widetilde{\alpha}_r}{\Sigma\widetilde{\alpha}_s}\simeq\matwo{\Sigma\widetilde{\alpha}_r}{0}\colon S^7\to P^4(3^r)\vee P^4(3^s) ~~\text{ for }r\leq s,\\
    &\mat{1_P}{0}{-i_4^s q_4^r}{1_P}\matwo{\Sigma\widetilde{\alpha}_r}{i_4^s\Sigma\alpha}\simeq\matwo{\Sigma\widetilde{\alpha}_r}{0}\colon S^7\to P^4(3^r)\vee P^5(3^s),
  \end{align*}
  where $i_4^s\colon S^4\to P^5(3^s)$ and $q_4^r \colon P^4(3^r)\to S^4$ are the canonical inclusion and pinch map, respectively.
Thus, after applying suitable elementary row matrix operations to the representation vector of the map $h$, we may assume that in the equation (\ref{eq:h}), $x_i=z_v=0$ for all $1\leq i\leq k$ and $1\leq v\leq n$, and there exists a unique nonzero $y_{u_0}$. Let $\bar{r}=r_{u_0}$ be the power index described in the statement of (\ref{PP=1:P4}). Since the homotopy cofibres of $\pm \Sigma \widetilde{\alpha}_{\bar{r}}$ are homotopy equivalent, we get the homotopy equivalence in (\ref{PP=1:P4}).

Secondly, by Lemma \ref{lem:Bockstein}, the assumption of (\ref{PP=1:S4}) implies that $y_u=z_v=0$ for any $1\leq u\leq m$ and $1\leq v\leq n$, and $x_i=\pm 1$ for some $1\leq i\leq k$. There is a homotopy 
\[\mat{1_4}{0}{-1_4}{1_4}\matwo{\Sigma \alpha}{\Sigma\alpha}\simeq \matwo{\Sigma\alpha}{0}\colon S^7\to S^4\vee S^4,\]
where $1_4$ denotes the identity on $S^4$.
Then after applying suitable elementary row matrix operations to the representation vector $M_h$, we may assume that 
\[x_1=1 \text{ and }x_2=\cdots=x_k=y_1=\cdots=y_m=z_1=\cdots=z_n=0.\] The statement (\ref{PP=1:S4}) is proved.

Thirdly, by Lemma \ref{lem:Bockstein}, the assumption of (\ref{PP=1:P5}) implies that $x_i=y_u=0$ for any $1\leq i\leq k$ and $1\leq u\leq m$, and  $z_v=\pm 1$ for at least one $v$. By the homotopy 
\[\mat{1_P}{0}{-B_4(\chi^r_s)}{1_P}\matwo{i_4^r(\Sigma \alpha)}{i_4^s(\Sigma \alpha)}\simeq \matwo{i_4^r(\Sigma\alpha)}{0}\colon S^7\to P^5(3^r)\vee P^5(3^s) \text{ for $r\geq s$},\]
we may assume that $z_v=1$ for exactly one $v=v_0$. Let $\bar{s}=s_{v_0}$ be the power index described in the statement of (\ref{PP=1:P5}), then we complete the proof of (\ref{PP=1:P5}).  
  \end{proof}

\end{proposition}

\begin{proof}[Proof of Theorem \ref{thm:main}]
  By Lemma \ref{lem:SM}, when localized away from $2$, $\Sigma M$ is homotopy equivalent to a wedge sum of spheres, Moore spaces and a copy of $C_h$.  Since $\PP$ acts trivially on the cohomology of spheres and Moore spaces, $\PP$ acts trivially on $H^3(M;\Z/3)$ if and only if $\PP$ acts trivially on $H^4(C_h;\Z/3)$. Combining Lemma \ref{lem:SM}, Propositions \ref{prop:PP=0} and \ref{prop:PP=1}, we then get the Theorem.
\end{proof}

We end this section by giving the $2$-local homotopy type of the double suspension $\Sigma^2 M_{(2)}$ in a special case.
\begin{theorem}\label{thm:SM-2-local}
	Let $M$ be a closed simply connected $7$-manifold with $H_\ast(M;\Z)$ given by (\ref{HMeq}), where $H$ and $T$ are both $2$-torsion-free. Suppose that $\Sq^2$ acts trivially on $H^5(M;\z{})$ and that $\Theta$ acts trivially on $H^\ast(M;\z{})$. 
  \begin{enumerate}
    \item There is a homotopy equivalence 
    \begin{equation}\label{SM5:2-local}
      \Sigma M_5\simeq_{(2)} \bigvee_{i=1}^{l-a}S^3\vee \bigvee_{j=1}^{k-b}S^4\vee \bigvee_{r=1}^{k-a}S^5\vee \bigvee_{s=1}^{l-b}S^6\vee \bigvee_{u=1}^a C^5_\eta\vee \bigvee_{v=1}^bC^6_\eta,
    \end{equation}
    where $0\leq a,b\leq \min\{k,l\}$. Moreover, $a=0$ if and only if $\Sq^2$ acts trivially on $H^2(M;\z{})$, and $b=0$ if and only if $\Sq^2$ acts trivially on $H^3(M;\z{})$.
    \item\label{SM:2-local} If $\Sq^2$ acts trivially on $H^2(M;\z{})$, then there is a homotopy equivalence 
    \[\Sigma^2 M\simeq_{(2)}\bigvee_{r=1}^k S^6\vee \bigvee_{s=1}^{l-b}S^7\vee \bigvee_{i=1}^{l-1}S^4\vee \bigvee_{j=1}^{k-b-1}S^5\vee \bigvee_{v=1}^{b-1}C^7_\eta\vee C_{\hbar},\]
where $\hbar\colon S^8\to S^4\vee S^5\vee C^7_\eta$ is given by the formal sum 
\begin{equation}\label{eq:hbar}
  \hbar\simeq_{(2)}x\cdot (\Sigma \nu')\eta+ y\cdot \Sigma^2\nu'+z\cdot i_5^\eta \Sigma^2\nu',
\end{equation}
    where $x\in\{0,1\}$, $(y,z)\in \{(0,0),(1,0),(2,0),(0,1)\}$.
  \end{enumerate}

	\begin{proof}
	(1) By assumption $H$ and $T$ do not contain $2$-torsion, implying that Moore spaces $P^m(H)$ and $P^n(T)$ become contractible after localization at $2$. The homotopy equivalence follows by Corollary \ref{cor:SM5-glob}.

    (2)  $\Sq^2$ acts trivially on $H^2(M;\z{})$ implying $a=0$ in (\ref{SM5:2-local}). 
  Consider the following compositions 
    \begin{align*}
      f_{S^5,r} &\colon S^7\xra{\Sigma f_5}W\twoheadrightarrow \bigvee_{r=1}^{k}S^5\twoheadrightarrow S^5_r,~~r=1,\cdots,k;\\
      f_{S^6,s} &\colon S^7\xra{\Sigma f_5}W\twoheadrightarrow \bigvee_{s=1}^{l-b}S^6\twoheadrightarrow  S^6_s,~~s=1,\cdots,l-b.
    \end{align*}
    Here $S^5_r$ and $S^6_s$ are the $r$-th $5$-sphere and $s$-th $6$-sphere in the wedge sums, respectively.
    Since $\Sq^2$ acts trivially on $H^6(\Sigma M;\z{})$ and $\Theta$ acts trivially on $H^\ast(M;\z{})$, the homotopy equivalence (\ref{SM5:2-local}) implies that $\Sq^2$ acts trivially on $H^6(C_{f_{S^6,s}};\z{})$ and $\Theta$ acts trivially on $H^5(C_{f_{S^5,r}};\z{})$, and hence the maps $f_{S^6,s}$ and $f_{S^5,r}$ are null homotopic for $s=1,\cdots,l-b$ and $r=1,\cdots,k$. Thus by Lemma \ref{STlemma}, there is a homotopy equivalence 
    \[\Sigma M\simeq_{(2)}\bigvee_{r=1}^{k}S^5\vee \bigvee_{s=1}^{l-b}S^6\vee C_{\varphi},\]
    where $C_\varphi$ is the homotopy cofibre of the  composition 
\[\varphi\colon S^7\xra{\Sigma f_5}\Sigma M_5\twoheadrightarrow \bigvee_{i=1}^{l}S^3\vee \bigvee_{j=1}^{k-b}S^4\vee \bigvee_{v=1}^bC^6_\eta.\] 
Consider the compositions
\[\varphi_{S^4,j}\colon S^7\xra{\varphi}\bigvee_{i=1}^{l}S^3\vee \bigvee_{j=1}^{k-b}S^4_j\vee \bigvee_{v=1}^bC^6_\eta\twoheadrightarrow S^4_i,\]
where $S^4_j$ is the $j$-th $4$-sphere of the wedge sum $\bigvee_{j=1}^{k-b}S^4_j$ for $1\leq j\leq k-b$.  
By similar arguments in the proof of Lemma \ref{lem:SM}, all cup products in $H^\ast(\Sigma M;\ZZ{2})$ are trivial implying that all cup products in $H^\ast(C_{\varphi_{S^4,j}};\ZZ{2})$ are trivial. Since $\pi_7(S^4)\cong_{(2)}\Z\lra{\nu_4}\oplus \Z/4\lra{\Sigma \nu'}$ and $\nu_4$ has the (2-local) Hopf invariant $1$, we apply \cite[Lemma 2.5]{CS} to get 
\[\varphi_{S^4,j}\simeq y_j\cdot \Sigma \nu', \text{ where $j=1,\cdots,k-b$}.\]
Similarly, we see that the composition
\[S^7\xra{\varphi} \bigvee_{i=1}^{l}S^3\vee \bigvee_{j=1}^{k-b}S^4\vee \bigvee_{v=1}^bC^6_\eta\twoheadrightarrow C^6_\eta\]
is generated by $i_4^\eta\Sigma \nu'\in\pi_7(C^6_\eta)$ in Lemma \ref{lem:pi7Ceta}. 
By \cite[Lemma 5.7, Proposition 5.8 and (5.9)]{Tod}, there are group isomorphisms 
\[\pi_7(S^3)\cong\z{}\lra{\nu'\eta_6}, \quad\pi_8(S^4)\cong \z{}\lra{(\Sigma \nu')\eta_7}\oplus\z{}\lra{\nu_4\eta_7}\] 
and homotopies $\nu'\eta_6\simeq \eta_3 \nu_4, \quad \eta_3(\Sigma \nu')\simeq \ast.$ Then by the Hilton-Milnor theorem, we may express the $2$-localized map $\varphi$ as the formal sum 
\[\varphi\simeq_{(2)}\sum_{i=1}^lx_i\cdot \nu'\eta+\sum_{j=1}^{k-b}y_j\cdot \Sigma\nu'+\sum_{v=1}^bz_v\cdot i_4^\eta\Sigma \nu'+\theta_\varphi,\] 
where $\theta_\varphi$ is a sum of Whitehead products,$x_i,z_v\in\z{}$ and $y_j\in\Z/4=\{0,\pm 1,2\}$ for each $i,j,v$. 
Applying the suspension functor to the above $2$-local formal sum of $\varphi$,  we get  
\[\Sigma \varphi\simeq_{(2)}\sum_{i=1}^lx_i\cdot (\Sigma \nu')\eta+\sum_{j=1}^{k-b}y_j\cdot \Sigma^2\nu'+\sum_{v=1}^bz_v\cdot i_5^\eta\Sigma^2 \nu',\]
where $x_i,z_v\in\z{}$ and $y_j\in\Z/4=\{0,\pm 1,2\}$ for each $i,j,v$.

There are homotopies
\begin{align*}
\mat{1_\eta}{0}{-1_\eta}{1_\eta}\matwo{i_5^\eta\Sigma^2\nu'}{i_5^\eta\Sigma^2\nu'}&\simeq \matwo{i_5^\eta\Sigma^2\nu'}{0},\\
  \mat{1_5}{0}{\mp 1_5}{1_5}\matwo{\Sigma^2\nu'}{\pm\Sigma^2\nu'}&\simeq\matwo{\Sigma^2\nu'}{0},\\
   \mat{1_5}{0}{-2}{1_5}\matwo{\Sigma^2\nu'}{2\Sigma^2\nu'}&\simeq\matwo{\Sigma^2\nu'}{0}.
\end{align*}
Here $1_n$ and $1_\eta$ denote the identity maps on $S^n$ and $C^{7}_\eta$, respectively.
Hence after applying a sequence of suitable elementary matrix operations to the representation vector 
\[M_{\Sigma \varphi}=[x_1\cdot (\Sigma \nu')\eta,\ldots,x_l\cdot (\Sigma \nu')\eta,\ldots,z_1\cdot i_5^\eta\Sigma^2\nu',\ldots,z_b\cdot i_5^\eta\Sigma^2\nu']^t,\]
we may assume that  
\begin{align*}
  &x_1\in \{0,1\},~~x_2=\cdots=x_l=0;\\
  &y_1\in\{0,1,2\}, ~~y_2=\cdots=y_{k-b}=0;\\
  &z_1\in\{0,1\}, z_2=\cdots=z_b=0.
\end{align*}
It follows that when localized at $2$, the suspension $\Sigma \varphi$ factors as the composition
\[\Sigma \varphi\colon S^8\xra{\hbar} S^4\vee S^5\vee C^7_\eta\hookrightarrow \bigvee_{i=1}^{l}S^3\vee \bigvee_{j=1}^{k-b}S^4\vee \bigvee_{v=1}^bC^6_\eta,\]
and by Lemma \ref{STlemma},  there is a homotopy equivalence 
\[\Sigma^2 M\simeq_{(2)}\bigvee_{r=1}^k S^6\vee \bigvee_{s=1}^{l-b}S^7\vee \bigvee_{i=1}^{l-1}S^4\vee \bigvee_{j=1}^{k-b-1}S^5\vee \bigvee_{v=1}^{b-1}C^7_\eta\vee C_{\hbar},\]
where $\hbar\colon S^8\to S^4\vee S^5\vee C^7_\eta$ is given by the formal sum 
\[\hbar\simeq_{(2)}x\cdot (\Sigma \nu')\eta+ y\cdot \Sigma^2\nu'+z\cdot i_5^\eta \Sigma^2\nu'\] 
for some $x,z\in\{0,1\}$ and $y\in \{0,1,2\}$. 
Consider the homotopies:
\begin{align*}
	& \mat{1_5}{0}{i_5}{1_\eta}\matwo{\Sigma ^2\nu'}{i_5^\eta\Sigma^2\nu'}\simeq\matwo{\Sigma^2 \nu'}{0}\colon S^8\to S^5\vee C^7_\eta,\\
	&\mat{1_5}{\overline{\zeta}}{0}{1_\eta}\matwo{2\Sigma^2 \nu'}{i_5^\eta\Sigma^2\nu'}\simeq\matwo{0}{i_5^\eta\Sigma\nu'}\colon S^8\to S^5\vee C^7_\eta,
\end{align*}
where the last homotopy is due to $\overline{\zeta}i_5^\eta \simeq 2\cdot 1_5$ by Lemma \ref{lem:piCeta}. 
Then we may assume that
\[x\in\{0,1\},\quad (y,z)\in \{(0,0),(1,0),(2,0),(0,1)\}. \]
Thus we get the $2$-local homotopy types of $C_\hbar$ and therefore that of $\Sigma^2M$ in (\ref{SM:2-local}). 
	\end{proof}
\end{theorem}

\section{Cohomotopy sets}\label{sec:chtp}
Let $M$ be a closed simply connected $7$-manifold with $H_\ast(M;\Z)$ given (\ref{HMeq}). 
This section is divided into three subsections to investigate the cohomotopy groups $\pi^k(M)$ for $k\geq 5$, $\pi^3(M)_{(\frac{1}{2})}$, and the $p$-local cohomotopy sets $\pi^4(M;\ZZ{p})=[M,S^4_{(p)}]$ for odd primes $p\geq 5$, respectively.

\subsection{Cohomotopy groups $\pi^k(M)$ for $k\geq 5$}
We shall determine the cohomotopy groups $\pi^k(M)$ in the stable range $k\geq 5$ in terms of the Postnikov tower of $S^k$ and the homotopy decomposition of $\Sigma M$. Recall that there are cohomotopy Hurewicz homomorphisms
\[h^k\colon \pi^k(M)\to H^k(M;\Z) \text{ for } k\geq 5.\]
First, we have the following result. 
\begin{lemma}\label{lem:chtp=6}
The following hold:
\begin{enumerate}
  \item $\pi^k(M)=0$ for $k=1$ or $k\geq 8$.
  \item\label{chtp=6} There are group isomorphisms \[\pi^7(M)\cong \Z\text{ and }\pi^6(M)\cong\z{1-\varepsilon},\]   
  where $\varepsilon=0$ if $\Sq^2$ acts trivially on $H^5(M;\z{})$, otherwise $\varepsilon=1$.
\end{enumerate}
 \begin{proof}
 (1) Since $\dim(M)=7$, we clearly have $\pi^k(M)=0$ for $k\geq 8$. Since $M$ is simply connected, we have $\pi^1(M)=H^1(M;\Z)=0$.

(2) As $M$ has dimension $7$, the cohomotopy Hurewicz  map (or the degree map) gives a group isomorphism $\pi^7(M)\cong H^7(M;\Z)\cong \Z$. 
 From \cite[Section 6.1]{Taylor12}, there is a short exact sequence of groups 
\[0\to H^7(M;\z{})/\Sq^2_\Z\big(H^5(M;\Z)\big)\to\pi^6(M)\xra{h^6}H^6(M;\Z)\to 0,\] 
where the integral Steenrod square $\Sq^2_\Z=\Sq^2\circ\rho_2$ with $\rho_2$ the mod $2$ reduction. Moreover, this exact sequence splits if and only if $\Sq^2_\Z$ and $\Sq^2$ have the same images in $H^7(M;\z{})$. 
Since $H^6(M;\Z)=0$, the splitting of the above short exact sequence is trivial, implying that 
\[\Sq^2_\Z\big(H^5(M;\Z)\big)=\Sq^2\big(H^5(M;\z{})\big), \quad \pi^6(M)\cong H^7(M;\z{})/\Sq^2\big(H^5(M;\z{})\big).\]
The group isomorphism $\pi^6(M)\cong\z{1-\varepsilon}$  in the Lemma follows immediately.
 \end{proof}
\end{lemma}

\begin{proposition}\label{prop:chtp=5:glob}
  Let $M$ be a smooth simply connected $7$-manifold with $H_\ast(M;\Z)$ given by (\ref{HMeq}). Set $\varepsilon=0$ if $M$ is spin and $\varepsilon=1$ otherwise. If $H$ is $2$-torsion-free, then there is a split exact sequence of groups
  \[0\to \z{1-\varepsilon}\to \pi^5(M)\xra{\jmath}\ker(\Sq^2_\Z)\to 0,\]
  where $\ker(\Sq^2_\Z)\cong H^5(M;\Z)$ and the composition 
  \[\pi^5(M)\xra{\jmath}\ker(\Sq^2_\Z)\xra{\cong} H^5(M;\Z)\]
  is the fifth cohomotopy Hurewicz homomorphism.
  \begin{proof}
    The splitting of the short exact sequence in the statement is straightforward because $\pi^5(M)$ is abelian and $H^5(M;\Z)\cong \Z^l\oplus H$ is $2$-torsion-free.
Consider the following two principal homotopy fibration sequences (the first two stages) in the Postnikov tower of $S^5$ (compare \cite[Page 122,Figure 1]{MT68}):
\begin{equation}\label{fib:2stages}
  \begin{aligned}
    &K(\Z/2,6)\xra{j_1}P_6S^5\xra{p_6}K(\Z,5)\xra{\Sq^2_\Z}K(\Z/2,7),\\
    & K(\z{},7)\xra{j_2} P_7S^5\xra{p_7}P_6S^5\xra{\overline{\Sq^2}}K(\Z/2,8),
  \end{aligned}
\end{equation}
where $P_kS^5$ is the $k$-th Postnikov section of $S^5$ such that the canonical map $S^5\to P_kS^5$ is $(k+1)$-connected for $k=6,7$, and $\overline{\Sq^2}\colon \Omega P_6S^5\to K(\z{},8)$ is the map such that $\Sq^2\simeq \overline{\Sq^2}\circ j_1$.
Since the above two stages (\ref{fib:2stages}) approximate $S^5$ in the stable range, we have the following induced exact sequences of groups:
  \begin{equation}\label{ES:chtp=5}
    \begin{aligned}
      &0\to [M,P_6S^5]\xra{(p_6)_\ast}H^5(M;\Z)\xra{\Sq^2_\Z}H^7(M;\z{}),\\
      &[M,\Omega P_6S^5]\xra{(\Omega\overline{\Sq^2})_\ast} H^7(M;\z{})\to [M,P_7S^5]\xra{(p_7)_\ast} [M,P_6S^5]\to 0
    \end{aligned}
  \end{equation}
The first exact sequence in (\ref{ES:chtp=5}) implies that 
  \[[M,P_6S^5]\cong \ker(\Sq^2_\Z)\cong H^5(M;\Z),\]
  where the last isomorphism holds because $H^5(M;\Z)$ is $2$-torsion-free.
  Since $M$ has dimension $7$, there is a group isomorphism $\pi^5(M)\cong [M,P_7S^5]$. Then from the second exact sequence in (\ref{ES:chtp=5}), we complete the proof of the Proposition by the \emph{claim} that the homomorphisms 
  \begin{align*}
    (\Omega\overline{\Sq^2})_\ast\colon [M,\Omega P_6S^5]\to H^7(M;\z{})\text{ and }\Sq^2\colon H^5(M;\z{})\to H^7(M;\z{})
  \end{align*}
  have the same cokernels $\z{1-\varepsilon}$.  
  
  If $M$ is nonspin, the Steenrod square  \[\Sq^2\colon H^5(M;\z{})\xra{(\Omega j_1)_\ast}[M,\Omega P_6S^5]\xra{(\Omega\overline{\Sq^2})_\ast}H^7(M;\z{})\]
  is an epimorphism implying that $(\Omega\overline{\Sq^2})_\ast$ is an epimorphism and hence has cokernel $0$.
 
  If $M$ is spin, we prove $\coker\big((\Omega\overline{\Sq^2})_\ast)\cong\z{}$ by showing that $\pi^5(M)\cong H^5(M)\oplus\z{}$. Let $\sk_3(M)$ be the $3$-skeleton of $M$. Since $M$ has dimension $7$, the cohomotopy sets $\pi^5(M)$ and $\pi^5(M/\sk_3(M))$ admit natural abelian group structure.  The cofibration $\sk_3(M)\hookrightarrow M\to M/\sk_3(M)$ implies group isomorphisms 
  \[\pi^5(M)\cong \pi^5(M/\sk_3(M))\cong \pi^6(\Sigma M/\Sigma \sk_3(M)).\]
   There is a homotopy cofibration 
   \[S^7\xra{\phi}\Sigma M_5/\Sigma \sk_3(M)\to \Sigma M/\sk_3(M),\] 
where $\phi$ is the attaching map that contains no Whitehead products.
   Collapsing cells of dimension at most $4$ in (\ref{SM5-glob}), we get a homotopy equivalence 
  \[\Sigma M_5/\Sigma \sk_3(M)\simeq \bigvee_{i=1}^{t+k} S^5\vee \bigvee_{j=1}^l S^6\vee P^6(H),\]
    where we assume that $P^5(T)$ becomes $\bigvee_{i=1}^t S^5$ after collapsing its $4$-skeleton. Since $\pi_7(P^6(p^r))=0 \text{ for }p\geq 3$, the space $P^6(H)$ retracts off $\Sigma M/\sk_3(M)$ by Lemma \ref{STlemma}.
    Since $M$ is smooth spin, $\Sq^2$ acts trivially on $H^5(M;\z{})$ and the secondary operation $\Theta$ that detects $\eta^2$ acts trivially on $H^\ast(M;\z{})$; consequently, the compositions of the form
    \[S^7\xra{\phi}\Sigma M_5/\Sigma \sk_3(M)\twoheadrightarrow S^6,\quad S^7\xra{\phi}\Sigma M_5/\Sigma \sk_3(M)\twoheadrightarrow S^5\]
    are both null homotopic, implying that 
    \[\Sigma M/\Sigma \sk_3(M)\simeq \bigvee_{i=1}^{t+k}S^5\vee \bigvee_{j=1}^l S^6\vee P^6(H)\vee S^8.\]
  Thus $ \pi^5(M)\cong \pi^6(\Sigma M/\Sigma \sk_3(M))\cong H^5(M;\Z)\oplus \z{}$ and we complete the proof of the claim.
\end{proof}
\end{proposition}

\begin{remark}\label{rmk:cohtp=5}
  Firstly, the exact sequence in Proposition \ref{prop:chtp=5:glob} not only implies the group structure of $\pi^5(M)$ but also explains it by the cohomotopy Hurewicz homomorphism. Secondly, the proof of Proposition \ref{prop:chtp=5:glob} shows that the cohomotopy group $\pi^k(M)$ in the stable range $\dim(M)\leq 2k-2$ only depends on the homotopy type of the quotient $\Sigma M/\Sigma\sk_{k-2}(M)$, a space typically far simpler than $\Sigma M$ itself.
\end{remark}

In the following two subsections, we shall investigate the cohomotopy set $\pi^k(M)$ with $k=3,4$ from the perspective of $p$-localization.

\subsection{The third cohomotopy group $\pi^3(M)$}\label{sec:chtp=3}
Consider the \emph{EHP} fibration sequence (compare \cite[Corollary 4.4.3]{N2})
\[\Omega^2 S^4\xra{\Omega H}\Omega^2 S^7\xra{P}S^3\xra{E}\Omega S^4\xra{H}\Omega S^7.\]
 Since $S^3$ is an $H$-space, $\pi^3(M)$ has a natural group structure; by \cite[Corollary 6.5]{HMRbook}, there are group isomorphisms
\[\pi^3(M)_{(\frac{1}{2})}\cong [M,S^3_{(\frac{1}{2})}]\cong [M_{(\frac{1}{2})},S^3_{(\frac{1}{2})}].\]
We shall combine the above \emph{EHP} fibration sequence and Theorem \ref{thm:main} to study the cohomotopy group $\pi^3(M)$.

\begin{proposition}\label{prop:chtp=3}
 The suspension $E_\ast\colon \pi^3(M)\to \pi^4(\Sigma M)$ is a group isomorphism in the following two cases:
 \begin{enumerate}
  \item the manifold $M$ is nonspin;
  \item all groups are localized away from $2$, that is, $E_\ast\colon \pi^3(M)_{\loca}\to \pi^4(\Sigma M)_{\loca}$ is an isomorphism.
 \end{enumerate}

\begin{proof}
Since $S^3$ is an $H$-space, the suspension map $E\colon S^3\to \Omega S^4$ has a left homotopy inverse; consequently, the suspension $E_\ast\colon \pi^3(M)\to \pi^4(\Sigma M)$ is injective. 

For the surjectivity, consider the EHP exact sequence of sets 
\begin{equation*}
	\pi^7(\Sigma^2 M)\xra{P_\ast}\pi^3(M)\xra{E_\ast}\pi^4(\Sigma M)\xra{H_\ast}\pi^7(\Sigma M).
\end{equation*}
Here we use the identification $[X,\Omega Y]=[\Sigma X,Y]$. 
By Lemma \ref{lem:chtp=6},  $\pi^7(\Sigma M)\cong \z{}$ for spin $M$ and $\pi^7(\Sigma M)=0$ for nonspin $M$.  Hence the suspension $E_\ast$ is surjective when localized away from $2$, or when $M$ is nonspin.

The obstruction for $E\colon S^3\to \Omega S^4$ to be an $H$-map is given by a map $S^3\wedge S^3\to \Omega S^4$ formed by $E$. Then Lemma \ref{lem:chtp=6} implies that any composition $M\to S^6\to \Omega S^4$ is trivial when $E_\ast$ is localized away from $2$, or if $M$ is nonspin. It follows that the suspension $E_\ast$ is a homomorphism, and hence an isomorphism in these two cases.
\end{proof}	 
\end{proposition}	

Recall that $X^8_r=P^5(3^r)\cup_{i_4\Sigma \alpha}e^8$ and $X^8(\widetilde{\alpha}_r)=P^4(3^r)\cup_{\Sigma \widetilde{\alpha}_r}e^8$ in Theorem \ref{thm:main}.
\begin{lemma}\label{lem:chtpgrps}
When localized away from $2$, there hold group isomorphisms 
\[[S^4\cup_{\Sigma\alpha}e^8,S^4]\cong \Z, \quad [X^8_r,S^4]=0,\quad [X^8(\widetilde{\alpha}_r),S^4]\cong \Z/3^{r-1}.\] 

	\begin{proof}
	Recall that $[P^5(p^r),S^4]\cong [P^6(p^r),S^4]=0$ for any odd prime $p$ by Lemma \ref{lem:Moore1}. 
When localized away from $2$, $\pi_8(S^4)=0$, we have the following exact sequences of groups:
		\begin{align*}
			&0\to [S^4\cup_{\Sigma\alpha}e^8,S^4]\to [S^4,S^4]\xra{(\Sigma\alpha)^\ast}[S^7,S^4],\quad 0\to [X^8_r,S^4]\to [P^5(3^r),S^4]=0,\\
			&0\to [X^8(\widetilde{\alpha}_r),S^4]\to [P^4(3^r),S^4]\xra{(\Sigma \widetilde{\alpha}_r)^\ast}[S^7,S^4]_{\loca}.
		\end{align*}
    Then we get the group isomorphisms in the Lemma.
	\end{proof}
\end{lemma}

Combining Theorem \ref{thm:main} and Lemma \ref{lem:chtpgrps}, we can compute the concrete group structure of $\pi^3(M)$ when localized away from $2$; the computational details are omitted.
\begin{corollary}\label{cor:chtp=3}
Let $M$ be given by Theorem \ref{thm:main}. If the condition of (\ref{SM:PP=1:P4}) holds, then there hold an isomorphism 
\[\pi^3(M)\cong_{\loca} \Z^k\oplus (H/(\Z/3^{\bar{r}}))\oplus \Z/3^{\bar{r}-1};\]
otherwise we have 
\[\pi^3(M)\cong_{\loca} \Z^k\oplus H\cong_{\loca} H^3(M;\Z).\] 

\end{corollary}

The following Example \ref{ex:notsurj} shows that the suspension $E_\ast\colon \pi^3(M)\to \pi^4(\Sigma M)$ is generally not surjective in the spin case. 

\begin{example}\label{ex:notsurj}
 The Poincar\'e duality complex $M=P^4(3)\cup_{i_3\alpha}e^7$ is a $S^3$-bundle over $S^7$, and hence a spin $7$-manifold by \cite[Theorem 5]{KS01}. Since $[P^{n+i}(3),S^n]=0$ for $n\geq 3$ and $i=1,2$ by Lemma \ref{lem:Moore1}, we have group isomorphisms 
 \begin{align*}
  \pi^3(M)&\cong \pi_7(S^3)\cong \z{},\\
  \pi^4(\Sigma M)&\cong \pi_8(S^4)\cong\z{}\oplus\z{},
 \end{align*} 
where the groups $\pi_7(S^3)$ and $\pi_8(S^4)$ refer to \cite[Proposition 5.8]{Tod}. Then it is clear that the suspension $E_\ast\colon \pi^3(M)\to \pi^4(\Sigma M)$ is not surjective.
\end{example}

\subsection{The fourth $p$-local cohomotopy set $\pi^4(M;\ZZ{p})$}\label{sec:chtp=4} 
Recall from  \cite[Chapter II, Theorems 5.1 and 5.3]{HMRbook} that the cohomotopy set $\pi^4(M)$ is the pullback of the canonical diagram
\[\pi^4(M;\Z_{\loca})\rightarrow \pi^4(M;\ZZ{0})\leftarrow \pi^4(M;\ZZ{2}),\]
and the induced map $\pi^4(M)\to \pi^4(M;\ZZ{2})$ and $\pi^4(M)\to \pi^4(M;\Z_{\loca})$ are injective. Hence $\pi^4(M)$ is closely related with $\pi^4(M;\ZZ{p})$. 
In this subsection, we shall apply the analysizing method in \cite[Section 6.3]{Taylor12} to study the \emph{$p$-local} cohomotopy set $\pi^4(M;\ZZ{p})=[M,S^4_{(p)}]$ for primes $p\geq 5$. 

Let $p\geq 5$ be a prime. From the homotopy groups  $\pi_i(S^4)$ for $i\leq 10$ computed by Toda \cite{Tod},  we have the following principal homotopy fibration sequence in the Postnikov tower of $S^4_{(p)}$ (compare \cite[Lemma 2.11]{GS21}):
\begin{equation}\label{fib:S4}
  K(\ZZ{p},3)\xra{\Omega \smallsmile^2}K(\ZZ{p},7)\xra{\partial}P_{10}S^4_{(p)}\xra{p_{10}} K(\ZZ{p},4)\xra{\smallsmile^2}K(\ZZ{p},8);
\end{equation}
where $P_{10}S^4_{(p)}$ is the $10$-th Postnikov section of $S^4_{(p)}$ such that the canonical map $\gamma_{10}\colon S^4_{(p)}\to P_{10}S^4_{(p)}$ is $11$-connected, and $\smallsmile^2$ is the cup squaring operation. In particular, for any CW-complex $X$ of dimension $\leq 10$, $(\gamma_{10})_\ast\colon \pi^4(X;\Z_{(p)})\to [X,P_{10}S^4_{(p)}]$ is a bijection and the composition 
\[\pi^4(X;\ZZ{p})\xra{(\gamma_{10})_\ast} [X,P_{10}S^4_{(p)}]\xra{(p_{10})_\ast}H^4(X;\ZZ{p})\]
is the $4$-th $p$-local cohomotopy Hurewicz map $h^4_{(p)}$.
Note that the canonical inclusion map $i_3\colon S^3_{(p)}\to K(\ZZ{p},3)$ is $10$-connected, implying that there is a bijection 
\[[K(\ZZ{p},3),K(\ZZ{p},7)]\to [S^3_{(p)},K(\ZZ{p},7)]\]
by \cite[Proposition 2.4.13]{Arkowitzbook}. It follows that the loop map $\Omega \smallsmile^2$ in (\ref{fib:S4}) is null homotopic, and hence the homotopy fibration sequence (\ref{fib:S4}) induces an exact sequence of sets  
\begin{equation}\label{ES:chtp=4}
  H^3(X;\ZZ{p})\xra{\ast}H^7(X;\ZZ{p})\to \pi^4(X;\ZZ{p})\xra{h^4_{(p)}}H^4(X;\ZZ{p})\xra{\smallsmile^2}H^8(X;\ZZ{p})
\end{equation}
for any CW-complex $X$ of dimension $\leq 10$.

Let $\mathcal{L}(-)=\mathrm{map}(S^1,-)$ be the free loop space functor on the category of topological spaces. Note that for any space $Y$, there is a homotopy fibration 
\[\Omega Y\to \mathcal{L}Y\xra{\mathrm{ev}_1}Y,\]
where $\mathrm{ev}_1(\omega)=\omega(1)$; moreover, the constant loop at a point $y\in Y$ defines a splitting section $s\colon Y\to \mathcal{L}Y$ such that $\mathrm{ev}_1\circ s\simeq 1_Y$. Hence there is an induced short exact sequence for each $i\geq 1$:
\[0\to \pi_i(\Omega Y)\to \pi_i(\mathcal{L}Y)\xra{(\mathrm{ev}_1)_\ast}\pi_i(Y)\to 0.\]

\begin{lemma}\label{lem:c_e}
 The homomorphism 
 \[ (\mathcal{L}\smallsmile^2)_\ast\colon H_7(\mathcal{L}K(\ZZ{p},4);\ZZ{p})\to H_7(\mathcal{L}K(\ZZ{p},8);\ZZ{p})\]
 is an isomorphism for each prime $p\geq 5$.
  \begin{proof}
  Consider the homotopy commutative diagram of homotopy fibrations induced by the canonical map $\gamma_{10}\colon S^4_{(p)}\to P_{10}S^4_{(p)}$:
  \[\begin{tikzcd}
  \Omega S^4_{(p)}\ar[r]\ar[d,"\Omega \gamma_{10}"]&\mathcal{L}S^4_{(p)}\ar[r,"\mathrm{ev}_1"]\ar[d,"\mathcal{L}\gamma_{10}"]&S^4_{(p)}\ar[d,"\gamma_{10}"]\\
  \Omega P_{10}S^4_{(p)}\ar[r]&\mathcal{L}P_{10}S^4_{(p)}\ar[r,"\mathrm{ev}_1"]&P_{10}S^4_{(p)}
  \end{tikzcd}\]
Then by the induced commutative diagram of short exact sequences of homotopy groups and the short Five-Lemma, we get that the induced homomorphism 
\[(\mathcal{L}\gamma_{10})_\ast\colon \pi_i(\mathcal{L}S^4_{(p)})\to \pi_i(\mathcal{L}P_{10}S^4_{(p)})\]
is an isomorphism for $i\leq 9$. Hence $\mathcal{L}\gamma_{10}\colon \mathcal{L}S^4_{(p)}\to \mathcal{L}P_{10}S^4_{(p)}$ is $8$-connected.
From \cite[the last line of page 20 and the table on page 21]{Ziller77}, we have 
\[H_6(\mathcal{L}S^4_{(p)};\ZZ{p})\cong \Z/2\otimes\ZZ{p}=0,\quad H_7(\mathcal{L}S^4_{(p)};\ZZ{p})=0,\]   
hence $H_6(\mathcal{L}P_{10}S^4_{(p)};\ZZ{p})=H_7(\mathcal{L}P_{10}S^4_{(p)};\ZZ{p})=0$.

By the Serre theorem (compare \cite[Theorem 6.4.2]{Arkowitzbook}) for the homotopy fibration $$\mathcal{L}P_{10}S^4_{(p)}\xra{\mathcal{L}p_{10}} \mathcal{L}K(\ZZ{p},4)\xra{\mathcal{L}\smallsmile^2} \mathcal{L}K(\ZZ{p},8),$$ 
the excision map $C_{\mathcal{L}p_{10}}\to \mathcal{L}K(\ZZ{p},8)$ is $10$-connected, where $C_{\mathcal{L}p_{10}}$ is the homotopy cofibre of $\mathcal{L}p_{10}$. Hence there is an exact sequence of homology groups 
\[H_7(\mathcal{L}P_{10}S^4_{(p)})\to H_7(\mathcal{L}K(\ZZ{p},4))\xra{(\mathcal{L}\smallsmile^2)_\ast} H_7(\mathcal{L}K(\ZZ{p},8))\to H_6(\mathcal{L}P_{10}S^4_{(p)}),\]
where the coefficients in $\ZZ{p}$ are omitted from the notation for brevity.
Since $H_6(\mathcal{L}P_{10}S^4_{(p)})=H_7(\mathcal{L}P_{10}S^4_{(p)})=0$, $(\mathcal{L}\smallsmile^2)_\ast$ is an isomorphism.
  \end{proof}
\end{lemma}

\begin{theorem}\label{thm:chtp=4}
  Let $X$ be a CW-complex of dimension $\leq 10$ and let $p\geq 5$ be an odd prime. 
  \begin{enumerate}[(1)]
    \item The cohomotopy Hurewicz map $h_{(p)}^4\colon \pi^4(X;\ZZ{p})\to H^4(X;\ZZ{p})$ is onto the subset of classes $u$ of $H^4(X;\ZZ{p})$ such that $u^2=u\smallsmile u=0$.
    \item For $u\in H^4(X;\ZZ{p})$ with $u^2=0$, let $e\in\pi^4(X;\ZZ{p})$ such that $h^4_{(p)}(e)=u$. Then there is a bijection between $(h_{(p)}^4)^{-1}(u)$ and the cokernel of 
    the homomorphism 
    \[\psi_e\colon H^3(X;\ZZ{p})\to H^7(X;\ZZ{p}),\quad \psi_e(v)=v\smallsmile e^\ast(\iota),\] 
    where $\iota\in H^4(S^4_{(p)};\Z_{(p)})$ is the fundamental class.
  \end{enumerate}
  \begin{proof}
    The first statement follows immediately from the exact sequence (\ref{ES:chtp=4}).

    For the second statement, since $(\Omega \smallsmile^2)_\ast$ in (\ref{ES:chtp=4}) is trivial, we apply \cite[Theorem 5.2 and Remark 5.3]{Taylor12} to the homotopy fibration sequence (\ref{fib:S4}) to get a bijection between $(h^{4}_{(p)})^{-1}(u)$ and the cokernel of the homomorphism 
    \[\phi_e\colon H^3(X;\ZZ{p})\to H^7(X;\ZZ{p}), \]
    where $\phi_e(v)$ is the composition 
    \[X\xra{\Delta}X\wedge X\xra{v\wedge e}K(\ZZ{p},3)\wedge S^4_{(p)}\xra{\mathfrak{D}'}K(\ZZ{p},7)\]
    for some map $\mathfrak{D}'$, where $u$ and $e$ are elements given in the statement. 
    Since $X$ has dimension $\leq 10$,  $v\colon X\to K(\Z_{(p)},3)$ factors through $i_3\colon S^3_{(p)}\to K(\Z_{(p)},3)$ and $\phi_e(v)$ factors as the composition 
\[X\xra{\Delta}X\wedge X\xra{v\wedge e}S^3_{(p)}\wedge S^4_{(p)}\xra{c_e}S^7_{(p)}\hookrightarrow K(\ZZ{p},7),\]
where $c_e$ is a map on $S^7_{(p)}$ of some degree. Here we abuse the notation $v$ to denote a cohomology class in $H^3(X;\ZZ{p})$ or a homotopy class in $\pi^3(X;\ZZ{p})$.  By the above factorization of $\phi_e(v)$, we have the formula
\[ \phi_e(v)=c_e\cdot (v\smallsmile e^\ast(\iota)),\]
where $\iota\in H^4(S^4_{(p)};\Z_{(p)})$ is the fundamental class. 

By \cite[Corollary 5.5]{Taylor12}, the map $\mathfrak{D}'$ in the factorization of $\phi_e(v)$, and hence the degree map $c_e$ is determined by the induced homomorphism $(\mathcal{L}\smallsmile^2)_\ast$ in Lemma \ref{lem:c_e}. Since $(\mathcal{L}\smallsmile^2)_\ast$ is an isomorphism,  the cokernel of $\phi_e$ is isomorphic to the cokernel the homomorphism $\psi_e$ given in the second statement. The proof of the Theorem is finished.
  \end{proof}
\end{theorem}

\begin{corollary}\label{cor:chtp=4}
  Let $M$ be a closed simply connected $7$-manifold and let $p\geq 5$ be a prime.
\begin{enumerate}[(1)]
  \item The cohomotopy Hurewicz map $h_{(p)}^4\colon\pi^4(M;\ZZ{p})\xra{} H^4(M;\ZZ{p})$ is surjective.
  \item For $u\in H^4(M;\ZZ{p})$, let $e\in \pi^4(M;\ZZ{p})$ be such that $h_{(p)}^4(e)=u$.  There is a bijection between $(h_{(p)}^4)^{-1}(u)$ and the cokernel of the homomorphism 
  \[\psi_e\colon H^3(M;\ZZ{p})\to H^7(M;\ZZ{p}),\quad \psi_e(v)=v\smallsmile e^\ast(\iota),\] 
  where $\iota\in H^4(S^4_{(p)};\ZZ{p})$ is the fundamental class. 
\end{enumerate}
\end{corollary}

\begin{proof}[Proof of Theorem \ref{thm:chtp}]
  The first two statements refer to Lemma \ref{lem:chtp=6}, the third statement refers to Proposition \ref{prop:chtp=5:glob}, the fourth statement refers to Proposition \ref{prop:chtp=3} and Corollary \ref{cor:chtp=3}, and the last statement is due to Corollary \ref{cor:chtp=4}.
  \end{proof}

\begin{Backmatter}
\bibliographystyle{amsplain}
\bibliography{7mflds}

\printaddress
\end{Backmatter}

\end{document}